\documentclass[a4paper,fleqn]{cas-sc}

\usepackage[authoryear,longnamesfirst]{natbib}

\def\tsc#1{\csdef{#1}{\textsc{\lowercase{#1}}\xspace}}
\tsc{WGM}
\tsc{QE}

\usepackage{tabularx}
\usepackage{pdflscape}

\begin{document}
\let\WriteBookmarks\relax
\def\floatpagepagefraction{1}
\def\textpagefraction{.001}

\shorttitle{\emph{In vitro} starch digestion kinetic}    

\shortauthors{F. Alonzo et~al.}  

\title [mode = title]{\emph{In vitro} starch digestion kinetic, formal analysis and new kinetic model}  

\tnotemark[1]

\tnotetext[1]{This research was co-funded by the Région Pays de la Loire (France).
This research benefits from funding from the French National Research Agency (ANR Glucoboost ANR-25-CE52-0594).}

\author[1]{Flavien Alonzo}[orcid=0000-0002-7972-3126]


\affiliation[1]{organization={INRAE UR BIA},
                addressline={Impasse Yvette Cauchois}, 
                city={Nantes},
                postcode={44300}, 
                state={},
                country={France}}

\author[1]{Yuzi Wang}
\fnmark[1]
\author[1]{Luc Saulnier}[orcid=0000-0002-0618-2672]

\author[2]{Jean-Philippe Ral}[orcid=0000-0002-3226-6204]


\affiliation[2]{organization={CSIRO Agriculture and food},
                addressline={GPO box 170}, 
                postcode={ACT 2601}, 
                postcodesep={}, 
                city={Canberra},
                country={Australia}}

\author[1]{Kamal Kansou}[orcid=0000-0002-6262-4447]
\cormark[1]
\ead{kamal.kansou@inrae.fr}

\cortext[cor1]{Corresponding author}
\fntext[fn1]{current professional address: John Wiley \& Sons, Shanghai, China}

\begin{abstract}
Describing mathematically the kinetics of the enzymatic hydrolysis of starch is a crucial part of 
nutrition and industrial applications studies. The majority of math models aim to fit the 
hydrolysis-extent curve ($X_t$), so-called digestogram, whose simple shape which can be fitted 
equally well by various math models. A relevant model of starch digestion kinetics should be more 
specific of the phenomenon. This study proposes a novel analytical framework that jointly fits the 
hydrolysis-extent curve ($X_t$) and the instantaneous reaction-rate 
coefficient ($k_a(t)$). 
A large dataset comprising 672 wheat-starch time-course 
measurements (including normal, waxy, and high-amylose genotypes), and additional 
experiments on other botanical starches was analysed. 
Standard First-Order and Weibull models were compared with two newly derived kinetic formulations 
\eqref{Eq.11} and \eqref{Eq.12} based on expected starch hydrolysis' phenomenon by which $k_a(t)$
tends with time towards an asymptotic-rate constant value, $k_e$, following a characteristic 
relaxation time ($\tau$). 
A mathematical proof of the existence of $k_e$, specific of the type of starch, is provided, based 
on an analysis of the Weibull's parameters.
Model ability to fit both $X_t$ and $k_a(t)$ was assessed by 
examining coefficients of determination ($R^2$) for both data types, and its adjusted version ($R^2_{adj}$) 
was used to compare performances between models. 
All models fit very well $X_t$ data, while $k_a(t)$ is much more discriminant of the model 
formulation. Wheat starch types are well distinguished, regarding ke, by the new models, with 
higher susceptibility to amylase for waxy samples ($\approx 5.7\times 10^{-3} min^{-1}$), than 
normal ($\approx 2.2\times 10^{-3} min^{-1}$) 
and high-amylose ($\approx 1.7\times 10^{-3} min^{-1}$) ones.
The new models outperformed the other models for monotonic evolution of $k_a(t)$, 
particularly digestion of high-amylose, and allow a more accurate estimation of $k_e$.
However, the analysis also highlighted limitations for modelling non-monotonic kinetics observed for 
waxy wheat starch digestion, suggesting a need for more sophisticated differential-equation 
descriptions for $k_a(t)$.
Overall, the combined fitting strategy provides a 
robust method for extracting physically meaningful kinetic parameters, facilitating 
deeper insight into the factors controlling starch digestion and supporting the 
development of more accurate predictive models for glycaemic response and food 
processing.
\end{abstract}


\begin{highlights}
\item Reaction-rate coefficient of starch hydrolysis by \(\alpha\)-amylase is often time-dependent
\item Reaction-rate coefficient tends towards a stable constant with time specific 
\item The combined fitting of hydrolysis extent and reaction-rate coefficient is proposed
\item New models more consistent with the experimental observations are proposed
\item New models outperformed usual models in the majority of cases
\end{highlights}

\begin{keywords}
 enzymatic hydrolysis \sep kinetic modeling \sep computational framework \sep amylase \sep first-order kinetic \sep wheat
\end{keywords}

\maketitle

\section{Introduction}

Starch is the natural product by which the plants store energy in the
form of polysaccharide for their metabolism and it is also a staple
dietary carbohydrate in human and animal diets. The enzymatic hydrolysis
of starch, which is called amylolysis, converts starch into glucose. It
is central in biological processes (human/animal digestion, grain
germination) and in industrial processes (malting, bioethanol
production). Determining the rate and the extent of amylolysis has
always been an important topic for the research area, especially for its
expected role in human nutrition.

The amylolysis is a heterogeneous hydrolytic reaction, that takes place
at the interface between the solid phase (starch) and the liquid phase
that conveys the enzyme. The main steps of the reaction include the
initial diffusion of amylase through the different levels of starch
organisation followed by the adsorption on the surface (formation of
enzyme-starch complex) and finally the hydrolysis of starch catalysed by
the enzyme \citet{Colonna1992LimitingFO}. The action of the amylolytic enzymes
is strongly affected by the starch structural features, which differ
greatly with the botanical source, macro-molecular composition or with
the transformation process (hydrothermal treatments, extrusion,
retrogradation) \citet{Zhang2015}, the native granular starch
\emph{in planta} displaying arguably a more complex hierarchical
structure than processed/cooked starch. Rapid and cheap \emph{in vitro}
assays mimicking starch digestion at different times have been developed
\citet{Zhang2015alt,Dhital2016}. A kinetics analysis of the
time-course measurements (sometimes called digestogram) can be performed
to identify the factors that affect the digestion process, this requires
the fitting of digestograms by a mathematical model \citet{Nguyen2018}. 
The first-order reaction exponential model adapted from \citet{Goi1997} and 
expressed in its integrated form \eqref{Eq.1} is widely reported in
the literature on this topic

\begin{equation}
    X_{t} = X_{\infty}\left( 1 - e^{- k t} \right),
    \label{Eq.1}
\end{equation}

In this model, $t$ denotes time, \(X_{t}\) the extent of hydrolysis
at time $t$, \(X_{\infty}\ \)the asymptotic value of \(X_{t}\) as
\(t \rightarrow + \infty\), and \(k\) the first-order reaction rate
constant. The model assumes that a starch fraction (\(X_{\infty})\) is
hydrolysed homogeneously by the enzyme; in other words, the reaction
proceeds at a constant rate relative to the remaining digestible starch.
This assumption is consistent with experimental results, in particular
the digestion of gelatinised purified starch, and widely adopted in 
the context of starch digestion analysis \citet{Butterworth2012,Edwards2014}.
Many studies have investigated the mechanisms underlying starch hydrolysis 
(e.g., \citet{Tahir2010,Warren2011,Baldwin2015,Li2019,Wang2022a}).
They showed that the hydrolysis of native starch or starch-containing matrices sometimes deviates 
from strict first-order reaction kinetics \citet{Butterworth2012,Edwards2014,Kansou2015a}.
These works consistently show that the hydrolysis of certain starches
follows an apparently biphasic pattern, providing a general picture of
amylase-driven starch digestion that can be divided into two main
stages:

\begin{itemize}
\item \textbf{Stage 1: Short term processes.} Diffusion and
adsorption of amylase on the substrate are governed by surface-dependent
factors.
\begin{enumerate}
\item[1a)] According to \citet{Warren2011}, the degree of order of the
glucan chains at the substrate surface limits amylase binding; once
amylase disrupts this order the limitation disappears.

\item[1b)] \citet{Baldwin2015} reported that, when immersed in water, solid
starch exposes surface glucan chains that are hydrolysed more rapidly
than the starch matrix itself. Damaged starch or starch with low degree
of order favours this rapid hydrolysis.

\item[1c)] Several factors can modulate the magnitude of mechanisms described in stages 1a and 1b,
the most important being the specific surface area of the substrate \citet{Wang2023}.\\
\end{enumerate}
\item \textbf{Stage 2: Longer term process.} As the short-term contributions
decline, starch is digested by amylase, at a rate that follows a pseudo
first-order kinetics law, i.e. a single ``inherent'' reaction rate
constant. This constant is mainly determined by the starch composition
and its internal organization as suggested in previous publications \citet{Zhang2015,Wang2023}.
\end{itemize}

Note that adverse phenomena capable of slowing or completely preventing
starch hydrolysis, such as amylose retrogradation, or enzyme inhibition,
are outside the scope of this work. Furthermore, the existence of an
intrinsic reaction-rate constant, as assumed in stage 2, has not
yet been demonstrated. More precisely, the factors governing this
constant remain unclear: it may reflect an intrinsic property of the
enzyme \citet{Edwards2014}, the substrate \citet{Wang2023} or
a combination of both. Consequently, an appropriate mathematical model
is required to capture the continuous transition from stage 1 to stage
2. Several models have been proposed in the literature to fit
accurately the digestograms, particularly variants of the exponential
model described in \eqref{Eq.1} (for a review see \citet{Nguyen2018}).
Unfortunately, the parameters of these models generally lack clear
physical meaning \citet{Nguyen2018}, particularly with respect to
the underlying postulated mechanism. Thus, many authors propose fitting
amylolysis curves with a double exponential model that represents the
sequential or parallel hydrolysis of two (or more) structural starch
fractions, each following its own first order kinetics 
\citet{Butterworth2012,Edwards2014,BelloPerez2019,Li2019}. 
This approach provides satisfactory fits to
experimental data in general, however, the physical basis for assuming distinct
starch fractions has been questioned \citet{Zhang2015}.
A more straightforward approach would be to fit an apparent reaction-rate 
'coefficient' over time, denoted $k_a(t)$. In a previous work, we
fitted starch digestion kinetics  using a Weibull
model \citet{Kansou2015b}, an extension of the exponential model \eqref{Eq.1} 
with an additional parameter for fitting $k_a(t)$ time-dependency while retaining 
simplicity and flexibility. For practical reason, the Weibull model is used in 
this work as a tool for analysing amylolysis kinetics. Nevertheless, the Weibull 
model presents drawbacks for fitting the reaction rate coefficient, as will be 
shown in this paper. The objective of this work is twofold:
(i) to provide a more reliable method for the analysis of starch digestion time-course 
measurement with fitting models, (ii) to propose a kinetic model for starch hydrolysis
that extends the simple first-order exponential formulation, aligning better
with the kinetics data and with the understanding of the starch-amylase system. 

To achieve this, we analysed a large dataset comprising hundreds of
\emph{in vitro} starch hydrolysis experiments with commercial $\alpha$-amylase,
encompassing native starches from cereals and tubers 
reported in multiple studies. Our proposed model offers a more dynamic perspective 
on $\alpha$‑amylase‑mediated starch degradation and should facilitate future 
investigations of this process.

\section{Materials and Methods}

The present study is based on a dataset obtained from \emph{in vitro}
starch hydrolysis assays performed for previous research works. In
particular the study of wheat starch hydrolysis kinetics can be found in
\citet{Wang2022b}.

\subsection{Materials}

\subsubsection{Native starch}

224 wheat starches (WS) from a 4-way MAGIC population \citet{Huang2012},
waxy and high amylose wheat starches (resp. WWS and HAWS) \citet{Ral2008} were
obtained from CSIRO (Australia). Potato starch, waxy maize starch
(Waxylis200) and high amylose maize (Eurylon7 (70\%), Eurylon5 (50\%))
were obtained from Roquette Frères (Lestrem, France), normal maize
starch (Meritena100) from Tereos (France), smooth pea starch (Nastar)
from Cosucra (Belgium), wrinkled pea starch (Amylose KG) from Stauderer,
(Germany), rice starch (Remy DR) from Remy industries (Belgium) and
amylose free potato starch from Lyckeby (Netherlands).



\subsubsection{Enzyme}

Porcine pancreatic $\alpha$-amylase (PPA) was used for the hydrolysis assays.
The enzymes were purchased from Neogen Corporation (Lansing, Michigan,
USA) for the wheat starch experiments, including waxy wheat and HAMY
wheat, and from Sigma (EU) for the remaining hydrolysis tests.

\subsection{Methods}

\subsubsection{Measurement of enzyme activity}

Enzyme solutions were prepared fresh before each hydrolysis and their
activity was checked using the Ceralpha Kit (Neogen Corporation,
Lansing, Michigan). PPA (2.5 mg/mL) was suspended in 20 mM sodium
phosphate buffer (pH 6.9) at 4℃ and stirred for 20 min. The solution was
then diluted 1:2000 and equilibrated for 20 min at 37℃. Thirty µL of the
diluted PPA solution were added to a pre-warmed 30 µL HR reagent
(blocked p-nitrophenyl $\alpha$-D-maltoheptaoside containing thermostable
$\alpha$-glucosidase) and incubated at 37℃ for 5 min. After this incubation, 50
µL reagent were transferred to 100 µL Ceralpha stopping buffer.
Absorbance was measured at 405 nm with a plate reader. Enzyme activity
was calculated using the following formula:

\begin{equation}
Act\left( \frac{U}{mL} \right) = \frac{{\mathrm{\Delta}OD}_{405}}{\varepsilon} \times \frac{1}{t} \times \frac{v}{V_{e}} \times F,\label{Eq.2}
\end{equation}

where \emph{$\Delta$OD\textsubscript{405}} is the absorbance difference between
sample and control, $\epsilon$ is the slope of the p-nitrophenol
calibration curve at 405 nm, \emph{V\textsubscript{e}} is the extraction
volume used for the assay, \emph{v} is the final assay volume, \emph{t}
is the reaction time (min) at which absorbance is read, and \emph{F} is
the dilution factor of the enzyme solution.

\subsubsection{Enzymatic hydrolysis of starch}

Starch (10 mg) and 1 mL of a 16 U/mL PPA solution were added separately
to the wells of a 96-well microplate and pre-warmed in a thermomixer
(Benchmark, USA) at 37℃ for 20 min. Then 1 mL PPA was added to each 10 mg
starch sample, the plate was sealed with a mat cover, and the mixture
was incubated at 37℃ with shaking at 1500 rpm. Aliquots (15 µL) were
withdrawn at 0, 20, 60, 120, 180, 240, 360, 1440 and 1800 min. Each
sample was originaly tested in triplicate or more. Note that in this work 
each replicate is fitted independently. 
This work exploits a dataset about wheat starch hydrolysis which
design of experiments can be found in \citet{Wang2022b}. This dataset includes
a total of 672 series of wheat starches (WS) hydrolysis time-course measurements (224
 distinct starch samples x 3 replicates), 14 replicates of Waxy Wheat Starch (WWS) and 14 of
High Amylose Wheat Starch (HAWS). Subsetting the database to the series without 
missing values reduces the dataset to 652 WS hydrolyses series. All replicated series
WWS and HAWS were selected. The starches from other biological sources have three replicates each.

\subsubsection{Reducing sugars measurements}

The extent of hydrolysis was monitored by measuring the appearance of
reducing sugars. Reducing sugars quantification was performed with the
3, 5-dinitrosalicylic acid (DNS) method \citet{Miller1959}. Aliquots (15 µL)
sampled at nine different time points were added to 100 µL of DNS
solution to stop the enzymatic reaction and to determine the reducing
sugar content. The mixture was then heated at 90°C for 10 min, and
allowed to cool at room temperature. Absorbance was measured at 540 nm
using a plate reader. Maltose was used as the calibration standard for
estimating the reducing sugar content.

\subsubsection{Extended procedure for fitting of starch degradation
models}

Starch degradation models are fitted with a least-squared algorithm that
simultaneously minimises two sets of residuals:
\begin{enumerate}
  \item \textbf{Hydrolysis‑extent residuals}, i.e., the difference between the
modelled extent of hydrolysis\,$X_t$ and the measured
value\,\({\widetilde{X}}_{i}\); 
\item \textbf{Instantaneous‑rate residuals},
i.e., the difference between the modelled instantaneous reaction
rate\,\(k_{a}(t)\) and the rate derived from the
data \({\widetilde{k}}_{a}(t_{i})\).
\end{enumerate}

The instantaneous rate\,\(k_{a}(t)\) is related to the hydrolysis extent
by \eqref{Eq.3a)}, and equivalently by \eqref{Eq.3b)}
\begin{align}
    X_{t} & = X_{\infty}\left( 1 - e^{- k_{a}(t)\cdot t} \right),\label{Eq.3a)}\tag{3a} \\
    \overset{t > 0}{\Leftrightarrow}k_{a}(t) & \overset{\mathrm{def}}{=} - \ln{\left( 1 - \frac{X_{t}}{X_{\infty}} \right)t^{- 1}}.\label{Eq.3b)}\tag{3b}
\end{align}

Thus, each measured hydrolysis extent \({\widetilde{X}}_{i}\), at time $t=t_i>0$, provides an indirect estimate of
the instantaneous rate:
\begin{equation}
  \widetilde{k}_{a}(t=t_{i}) = \ln{\left( 1 - \frac{\widetilde{X}_i}{X_{\infty}} \right)t_i^{- 1}}.\label{Eq.3c}\tag{3c}
\end{equation}
\setcounter{equation}{3}

Substituting \(X_{t}\) in \eqref{Eq.3b)} by the model expression gives a rate
expression that depends on the model parameters. In this work, we consider any kinetic model.
Thus, the first-order model \eqref{Eq.1} obviously yields $k_a^{first-order}(t)=k$. 
The Weibull model \eqref{Eq.4} provides a better example. It introduces two extra parameters
\(k_{w}\) and \(h\) while retaining the same notations as in \eqref{Eq.1}.

\begin{equation}
    X_{t} = X_{\infty}\left( 1 - e^{- k_{w}t^{1 - h}} \right),\ \ with\ h < 1.\label{Eq.4}
\end{equation}

Note that \(0 < h < 1\) gives a classic exponential saturation curve; while \(h < 0\) 
gives a S-shaped curve. Substituting the expression for $X_t$ from \eqref{Eq.3b)} for the 
Weibull model yields 
\begin{equation*}
    k_{a}^{Weibull}(t) = k_{w}t^{- h}.
\end{equation*}

This makes the modelling time-dependency of $k_a(t)$ explicit. In Appendix.A, a list of model 
equations used in the paper and in the literature is provided with their corresponding
  $k_a^{model}(t)$.

\subsubsection{Overall fitting criterion}

The fitting procedure minimises the following least-squared objective:
\begin{equation*}
    minimize\ \sum_{}^{}\left( X^{model}\left( t_{i} \right) - {\widetilde{X}}_{i} \right)^{2} + \lambda\sum_{}^{}\left( k_{a}^{model}\left( t_{i} \right) - {\widetilde{k}}_{a}\left( t_{i} \right) \right)^{2}.
\end{equation*}

\(where\ \lambda\) is a numerical hyperparameter (equal to 1000 in this
work) that scales the rate residuals so that they are comparable in
magnitude to the extent residuals.

Model performance is assessed with the coefficient of determination
(\(R^{2}\)) calculated separately for the hydrolysis extent and the
instantaneous rate fits to verify that both curves are well reproduced.
The adjusted coefficient of determination (\(R_{adj}^{2}\)) is also 
reported to account for the number of parameters of the various models.

\section{Results}

\subsection{Analysing the instantaneous reaction rate ($k_a$) time-dependency}

\begin{figure}
\centering
\includegraphics[width=\linewidth]{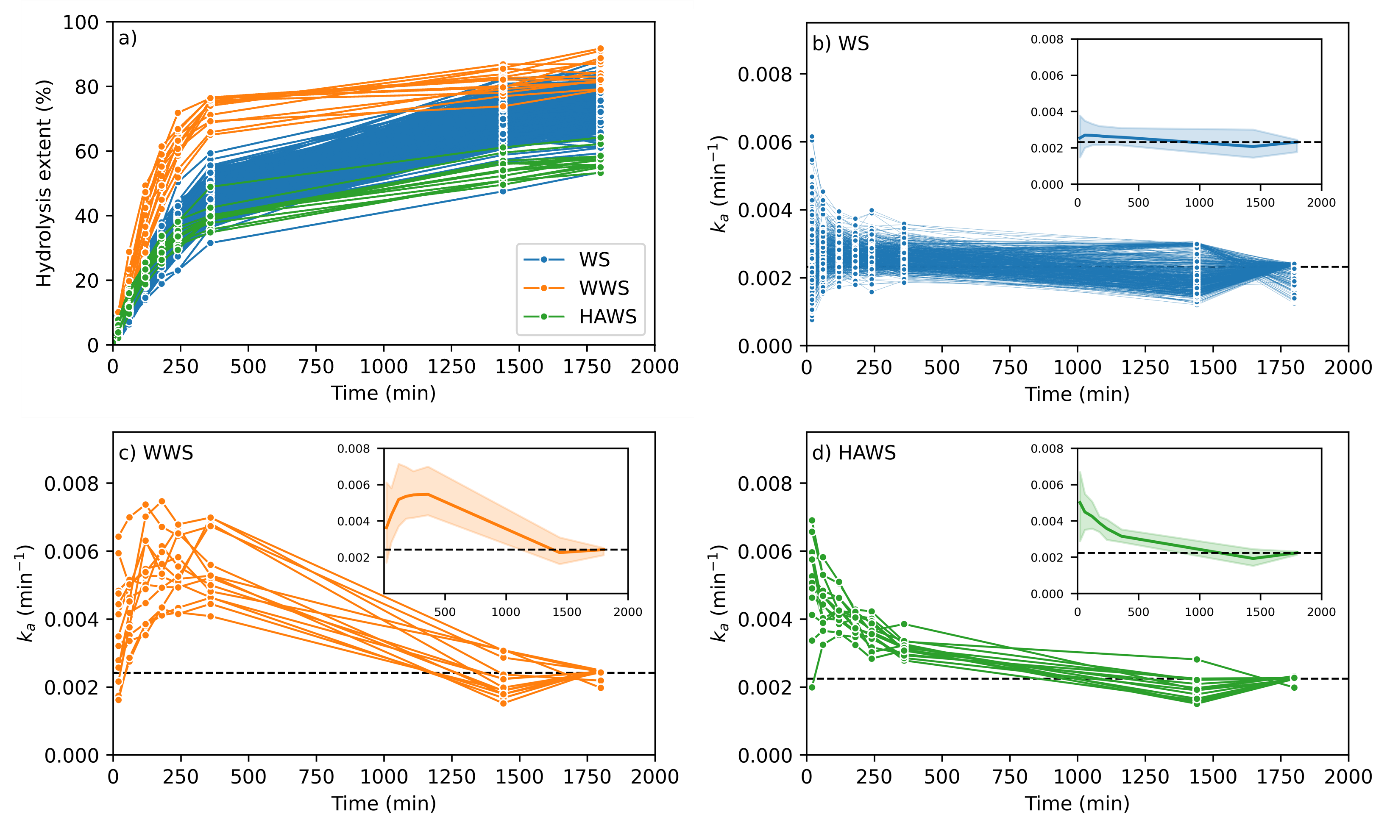}
\caption{Starch hydrolysis by PPA time-course measurements for
wheat starches (data from \citet{Wang2022b}). Panel a: hydrolysis extent for 680 wheat starches split into 652 normal
starches (WS), 14 waxy starches (WWS), and 14 high amylose starches
(HAWS, 72.9\% amylose).
Panels b-d) instantaneous reaction rate coefficient $k_a$
derived from the data in panel a, shown separately for WS (b), WWS (c),
and HAWS (d).
In each subplot the individual curves are plotted in the main graph; the
horizontal dashed line indicates the apparent
steady-state of $k_a(t)$ at long time and is provided as a
visual guide. In the small graphs, the thick line represents the mean curve and the shaded area shows the 90\%
percentile interval.}
\label{Fig.1}
\end{figure}

In the general case, starch is rapidly hydrolysed by amylase during an
initial phase, after which the reaction slows down progressively as the
substrate becomes depleted, eventually reaching a negligible rate. This
behaviour is illustrated in Fig.\ref{Fig.1}a which shows amylolysis kinetics
curves for three wheat‑starch types that were measured in a previous
work \citet{Wang2022b}: normal wheat starch (WS), waxy wheat starch (WWS), and
high amylose wheat starch (HAWS). The dataset comprises 652 time-course
measurements for WS and 14 measurements each for WWS and HAWS. As
expected, WWS, is almost completely hydrolysed by t=1800 min, whereas
HAWS, retains a larger amount of residual starch at the same time point;
WS displays intermediate behaviour ($\approx$\,70\,\% hydrolysis). The
variability of the WS samples digestion reflects the diversity of
wheat genotypes and growing conditions \citet{Wang2022b}, 
whereas the HAWS and WWS measurements are replicates
obtained over a three-month period. WS samples also cover a broad range of
granule-size distributions, while the ranges of compositional and
structural ranges are narrower in comparison \citet{Wang2022b}. 
Measurements of the corresponding instantaneous reaction-rate 
coefficient $\widetilde{k}_a(t=t_i)$ are plotted in Fig.\ref{Fig.1}b-d 
for WS, WWS, and HAWS, respectively. These curves were obtained 
by applying \eqref{Eq.3c} (see Methods section).

Examining the estimated instantaneous reaction rate\,\(k_{a}(t)\) over the
course of hydrolysis for the three types of wheat starch reveals three
salient monotonic kinetics (Fig.\ref{Fig.1}b-d).

\begin{enumerate}
\item
  \textbf{Rapid decline:} the initial value
  of\,\(k_{a}(t)\) is high, then drops sharply with time and stabilises at
  a non‑zero plateau. This pattern is predominant for
  HAWS digestion (Fig.\ref{Fig.1}d), and for many WS samples.
\item
  \textbf{Nearly constant:} \(k_{a}(t)\) remains
  roughly constant throughout the reaction, which is consistent with 
  the first‑order reaction kinetics and observed for many WS samples.
\item
  \textbf{Slow rise:} in a few cases for WS (Fig.\ref{Fig.1}b), the 
  initial\,\(k_{a}\) is low and gradually increases until
  it reaches the same range as the other measurements. This yields a 
  S‑shaped digestion curve less frequently reported in the literature; it
  typically occurs when the starch is not directly accessible to amylase, 
  for example when starch granules are entrapped within cell‑wall matrices
  \citet{PallaresPallares2018}.
\end{enumerate}

A different dynamic is illustrated by WWS digestion (Fig.\ref{Fig.1}c).
Here the evolution of \(k_{a}\) is non‑monotonic: the rate rises
quickly during the first few hours of the assay, then declines
during the rest of the hydrolysis.

\begin{figure}
\includegraphics[width=\linewidth]{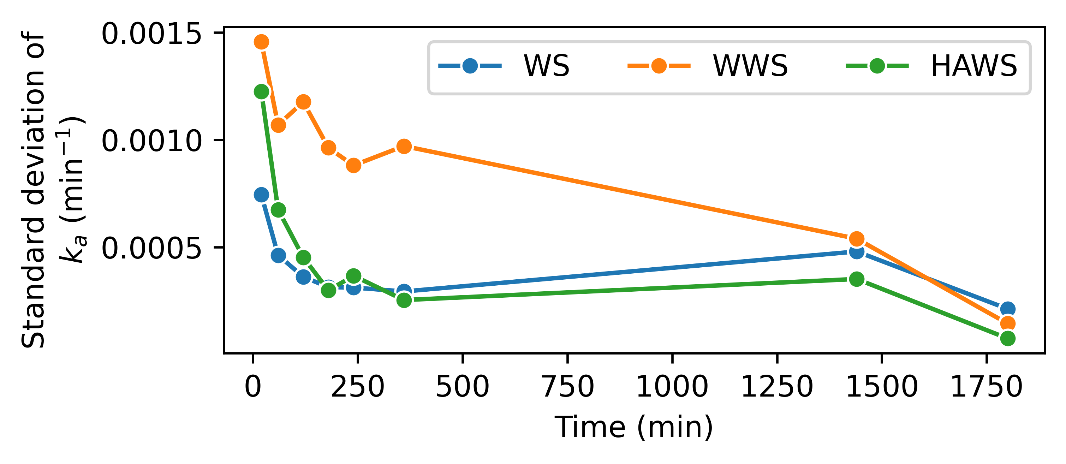}
\caption{Standard deviations of the instantaneous reaction rate
coefficient (\(k_{a}\)) against time, for WS (wheat starch), WWS (Waxy
wheat starch), and HAWS (high-amylose wheat starch) samples.}
\label{Fig.2}
\end{figure}

In Fig.\ref{Fig.1}b-d the variability of $k_a$ values changes
markedly over time. To confirm this visual impression, the standard
deviations of \({\widetilde{k}}_{a}\) was plotted against the time in
Fig.\ref{Fig.2}. This plot clearly shows a decrease in the standard
deviations as the hydrolysis progresses for all wheat starch
samples.

It is commonly put forward in publications that starch digestion by 
amylase proceeds in two
apparent phases. In the rapid‑digestion phase, ``dangling'' glucan
chains exposed at the substrate surface facilitate enzyme access,
whereas the slower phase is limited by the reduced surface area
available for enzyme binding. These surface‑dependent mechanisms
dominate the early stage of hydrolysis and are gradually overtaken by
the breakdown of the more compact starch matrix. This is consistent 
with our results for many wheat starch and for high‑amylose wheat 
starch for which the instantaneous reaction rate is first highly 
variable and tends to decrease until it converges over time to 
a constant value that follows first‑order kinetics (Fig.\ref{Fig.1}b,d and Fig.\ref{Fig.2}); 
we refer to this value as the ''inherent'' reaction‑rate constant, denoted\,\(k_{e}\).
The decreasing variability of \(k_{e}\) is also consistent with this view.
For waxy wheat starch (WWS) the dynamics are more complex because the
rate is non‑monotonic, yet a convergence of the instantaneous reaction rate
toward a constant is very likely.

\subsection{Revealing the ''inherent'' reaction rate constant \(k_{e}\)} 
\label{revealing ke}

The results presented above, together with findings from other studies,
show that the instantaneous reaction rate\,\(k_{a}(t)\) tends to
converge toward a constant value for various type of wheat starch 
(WS, WWS, HAWS). To our knowledge, the existence of such an “inherent” reaction 
rate constant has never been formally demonstrated in the literature. In the 
present work we provide an indirect but formal proof by fitting the starch hydrolysis 
curves with the Weibull model (\eqref{Eq.4}; see Appendix.B). 

The Weibull model \eqref{Eq.4} contains three parameters: the total
hydrolysis extent\,\(X_{\infty}\), a reaction‑rate constant\,\(k_{w}\),
and a shape parameter\,\(h\). In principle these parameters are
independent, but when the Weibull model is fitted to digestograms a
strong correlation between\,\(k_{w}\) and\,\(h\) might be observed, as previously
reported in \citet{Olawoye2020,Wang2022b}. This correlation can be expressed as
\begin{equation}
    h = a \times \ln\left( k_{w} \right) + b,\ \ a,b\ constants.\label{Eq.5}
\end{equation}

\begin{figure}
\centering
\includegraphics[width=\linewidth]{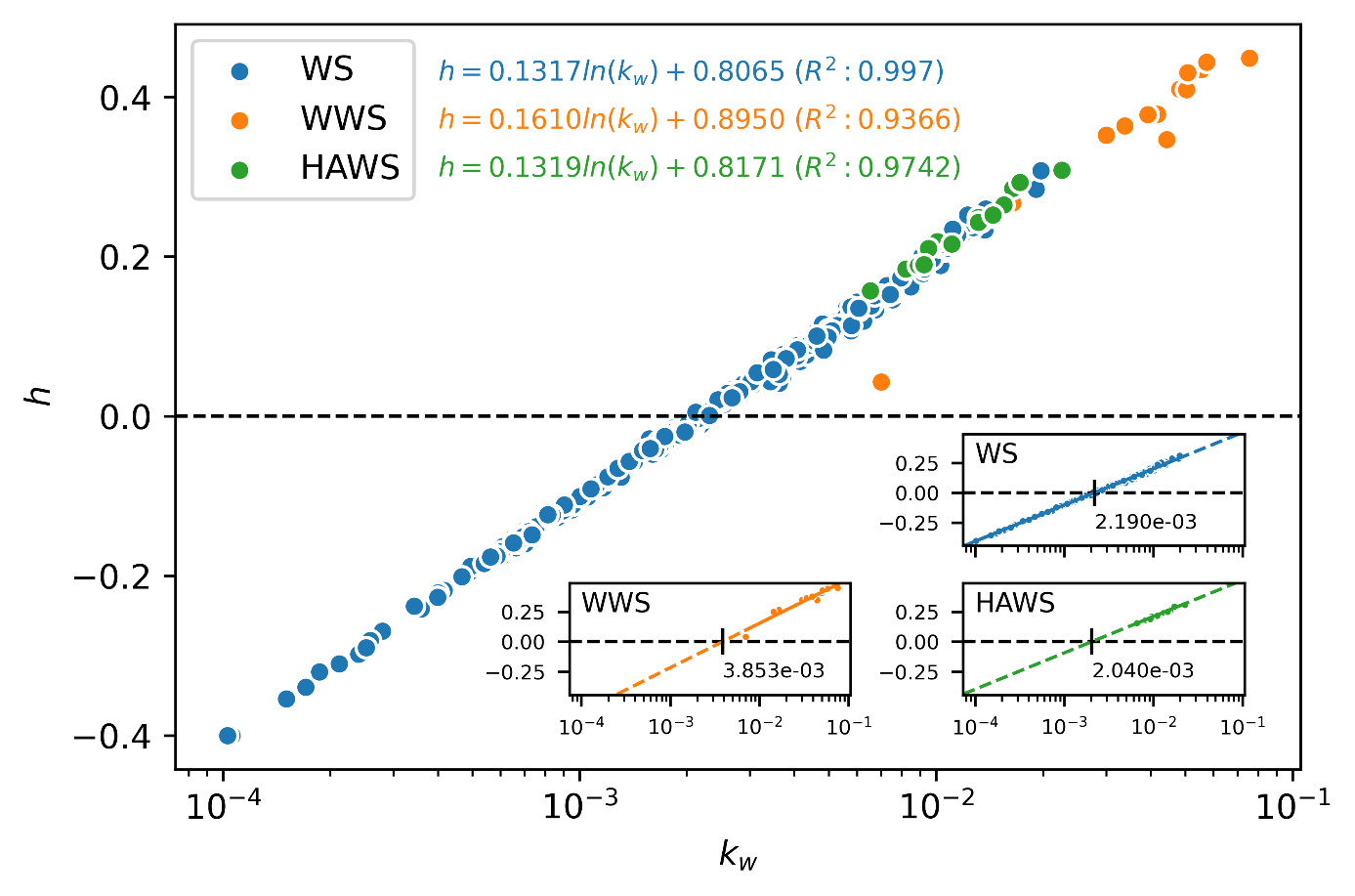}
\caption{Relationship between the parameters $k_w$ and $h$ of the Weibull model \eqref{Eq.4}
 fitted by minimizing \eqref{Eq.7}. Parameters’ distribution related to WS samples 
 are in blue, WWS in orange, and HAWS in green. Explicit relationships between $h$ and $ln(k_w)$ 
 are given at the top of the plot with their $R^2$-scores. Smaller plots show the curve 
 associated with each sample type and provide $k_w^0$, i.e., the estimated value of the inherent 
 reaction rate constant $k_e$.}
\label{Fig.3}
\end{figure}

Fig.\ref{Fig.3} illustrates this relationship for the WS, WWS and HAWS samples, showing
an excellent fit of the data to \eqref{Eq.5}. In \eqref{Eq.5} the particular value
of\,\(k_{w}\) for which\,\(h = 0\) is of special interest because, in
the Weibull model,\,\(h = 0\) reduces the equation to the first‑order
model \eqref{Eq.1} with the rate constant\,\(k = k_{w}\). We denote this
value by\,\(k_{w}^{0}\), which can be written as
\(k_{w}^{0} = \exp\left( - \frac{b}{a} \right).
\) According to Fig.\ref{Fig.3}, the estimated\,\(k_{w}^{0}\) values are:

\begin{itemize}
\item
  \textbf{WS:} \(k_{w}^{0} = 2.190 \times 10^{-3}\)\,(\({min}^{-1}\));
  95\,\%\,CI\,{[}2.073\,;\,2.322{]}\,×\,10\({}^{-3}\),
\item
  \textbf{WWS:} \(k_{w}^{0} = 3.853\times 10^{-3}\)\,(\({min}^{-1}\)); 95\,\%\,CI\,{[}\(0.671\); \(13.496\){]}
  \(\times10^{-3}\),
\item
  \textbf{HAWS:} \(k_{w}^{0} = 2.040\times 10^{-3}\)\,(\({min}^{-1}\)); 95\,\%\,CI\,{[}\(0.608\); \(5.463\){]}
  \(\times10^{-3}\).
\end{itemize}

It can be shown that, for a series of digestion kinetics of a given starch, the existence of a
finite limit of the instantaneous reaction‑rate
coefficient\,\(k_{a}(t)\) as\,\(t \rightarrow \infty\); i.e., \(k_{e}\) is mathematically
equivalent to the relationship expressed in \eqref{Eq.5}. The complete
proof is provided in Appendix.B. This result can be summarized as
\begin{equation}
    \lim_{t \rightarrow + \infty}{k_{a}(t)} = k_{e}\ constant\ \overset{}{\Leftrightarrow}\ \exists(a,b)\ constants:\ h = a \times \ln\left( k_{w} \right) + b.\label{Eq.6}
\end{equation}

\subsubsection{Principle of the demonstration}
\eqref{Eq.6} follows from fitting the hydrolysis‑extent curves for a given starch 
with the Weibull model \eqref{Eq.4}. The fitting procedure minimises the residuals between the
modelled extent\,\(X(t=t_{i})\) and the measured
values\,\({\widetilde{X}}_{i}\) at each sampling time\,\(t_{i}\).
Formally, the calibration problem is minimizing the function:
\begin{equation}
    J\left( k_{w},h \right) = \ \frac{1}{2}\sum_{i}^{}\left( k_{w}t_{i}^{1 - h} - {\widetilde{y}}_{i} \right)^{2},\label{Eq.7}
\end{equation}

where
\({\widetilde{y}}_{i} = \  - \ln\left( 1 - \frac{\ {\widetilde{X}}_{i}}{X_{\infty}} \right) = {\widetilde{k}}_{a}(t_{i}) \times t_{i}\).
Because the Weibull model gives increasing weight to later points, the
residual associated with the latest measurement (\(t_{\max}\)) dominates
the optimisation, enforcing
\(k_{w}t_{\max}^{1 - h} - {\widetilde{y}}_{\max} \approx 0\). Using
\eqref{Eq.5}, this condition can be rewritten as
\begin{equation}
    h = {\frac{1}{\ln\left( t_{\max} \right)}\ln\left( k_{w} \right) - \frac{\ln\left( {\widetilde{k}}_{a}\left( t = t_{\max} \right) \right)}{\ln\left( t_{\max} \right)}},\label{Eq.8}
\end{equation}

Since the experimental protocol uses a common final measurement time
(\(t_{\max} = 1800\text{\:\,min}\)) for all samples, \eqref{Eq.8} implies that
\({\widetilde{k}}_{a}(t_{\max})\) must approach a constant value;
otherwise, the intercept\,\(b\) in \eqref{Eq.5} would not be constant.
Consequently, a finite limit for \(k_{a}(t)\) is a necessary and
sufficient condition for the validity of \eqref{Eq.5}.

This limit value is denoted by \(k_{e}\):
\begin{equation*}
    k_{e} = \underset{t \rightarrow \infty}{\lim}k_{a}(t).
\end{equation*}

When \(h = 0\) the Weibull model reduces to the first‑order
(exponential) model, and the corresponding value of the Weibull rate
constant is 
\begin{equation}
    \widetilde{k_e} = k_{w}(h = 0) = k_{w}^{0} = e^{\frac{- b}{a}}.\label{Eq.9}
\end{equation}

Thus, an estimate of \(k_{e}\) can be obtained from the linear relation
between \(h\) and \(\ln(k_{w})\) \eqref{Eq.5}. Note that the confidence
intervals for the estimate derived from \eqref{Eq.9} decrease only
as\,\(1/\sqrt{n}\), where\,\(n\) is the number of time‑course
measurements; the detailed derivation of these confidence intervals is given in
Appendix.C.

\subsection{Determination of the inherent reaction‑rate constant\,\(k_{e}\) for various native starches}

\begin{figure}
\centering
\includegraphics[width=\linewidth]{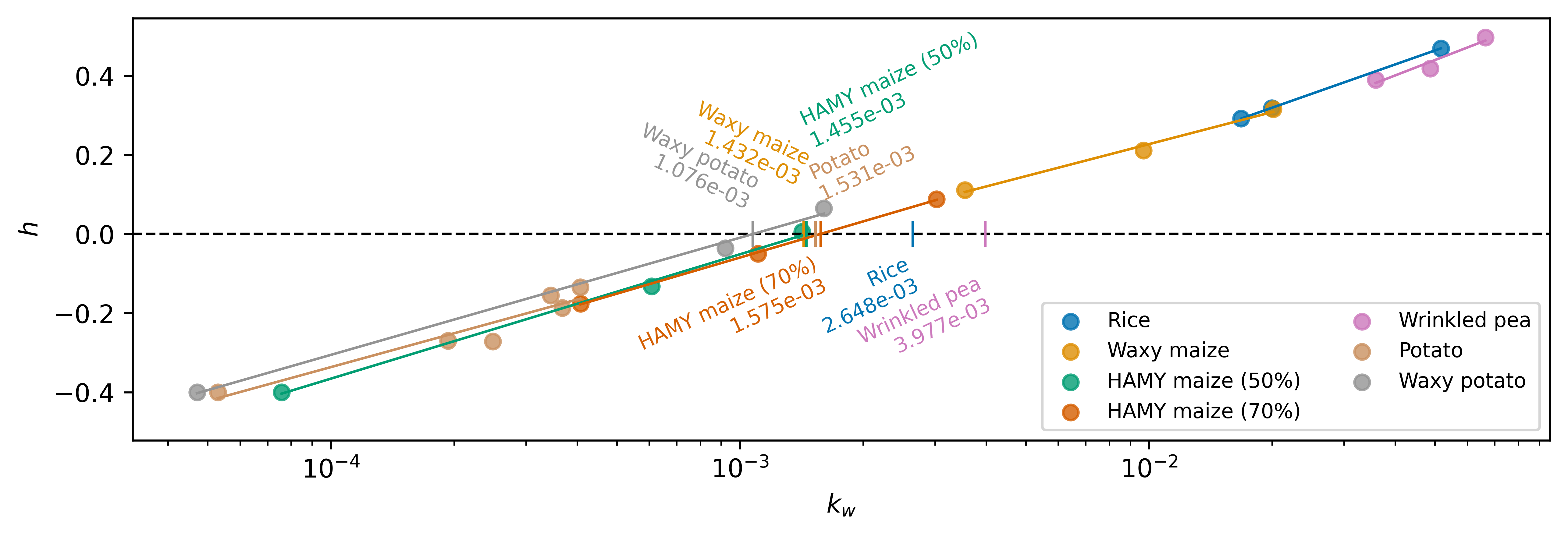}
\caption{Linear relationship between the Weibull
parameters\, $h$ and $ln(k_{w})$ for starches from different botanical origins. 
The corresponding \(k_{w}^{0}\) value, i.e., estimation of the inherent reaction‑rate
constant \(k_{e}\), is reported alongside each starch type.}
\label{Fig.4}
\end{figure}

Starches from rice, wrinkled pea, normal and waxy potato, and waxy and
high‑amylose (HAMY) maize were hydrolysed with amylase in three
independent experiments each. Each time‑course curve was fitted with the
Weibull model \eqref{Eq.4}, which yields the parameters\,\(k_{w}\)
and\,\(h\). Plotting\,\(h\) against\,\(\ln(k_{w})\) (Fig.\ref{Fig.4}) reproduces
the linear relationship described in \eqref{Eq.5}, allowing the computation 
of the "inherent" rate constant, \(k_{e}\) by means of \eqref{Eq.9}.

The estimated\,\(k_{e}\) values reveal clear differences among starch
types. Wrinkled pea and WWS exhibit the highest \(k_{e}\)
(\(0.00398\,{min}^{- 1}\), and \(0.00385\,{min}^{- 1}\) respectively), whereas 
waxy potato and native potato starches display the lowest\,\(k_{e}\)
(\(0.00108\,{min}^{- 1}\), and \(0.00114\,{min}^{- 1}\) respectively). 
In general, starches from cereals and legumes (rice, wheat, wrinkled pea) 
have higher\,\(k_{e}\) values than those from high‑amylose
plants and tubers (waxy and HAMY maize, potato, waxy potato) (Fig.\ref{Fig.4}).
This ranking agrees with the expected susceptibility of these starches to amylase, 
for example waxy wheat starch is readily hydrolysed, while high‑amylose and tuber 
starches are more resistant to digestion.

These results suggest that\,\(k_{e}\) captures an intrinsic and stable 
property of the starch-amylase system and should therefore be included 
as a model parameter.
Consequently, it appears beneficial to consider both the hydrolysis 
extent\,\(X_t\) and the instantaneous reaction‑rate coefficient\,
\(k_{a}(t)\), when fitting model starch hydrolysis.

\subsection{New models using the amylolysis ``inherent'' reaction rate
constant ($k_e$)}

The analysis of amylolysis kinetics described above suggests the
following requirements for semi‑empirical models of starch digestion by
amylase:

\begin{itemize}
\item
  Explicitly describe the continuous evolution of the reaction‑rate
  coefficient over time.
\item
  Honour the initial
  conditions\,\(X_t(0) = 0\)\,and\,\(k_{a}(0) = k_{\text{initial}}\).
\item
  Converge to a finite, constant value
  as \(t \rightarrow \infty\)\,(i.e.,\,\({\lim}_{t \rightarrow \infty}X_t = X_{\infty}\)
  and \({\lim}_{t \rightarrow \infty}k_{a}(t) = k_{e}\)).
\item
  Use the smallest possible number of parameters.
\item
  Reduce to the classic first‑order exponential model when appropriate
  (i.e., be a generalisation of that model).
\item
  Remain flexible enough to capture both rapid‑to‑slow and slow‑to‑rapid
  (S‑shaped) kinetic patterns.
\end{itemize}

\subsubsection{New model with the required specificities}

The model builds on the conventional hydrolysis dynamics \eqref{Eq.10a} and
couples it to a linear first‑order differential equation that governs
the evolution of $k_t$, the instantaneous reaction‑rate coefficient \eqref{Eq.10b}:
\begin{align}
    \dot{X_{t}} & = k(t) \times \left( X_{\infty} - X_{t} \right),\label{Eq.10a}\tag{10a} \\
    \dot{k_{t}} & = \frac{1}{\tau} \times \left( k_{e} - k_{t} \right).\label{Eq.10b}\tag{10b}
\end{align}

\begin{itemize}
\item
  \textbf{Initial conditions:} \(X_t(0) = 0\) and \(k_t(0) = k_{0}\).
\item
  \textbf{Parameters:} \(X_{\infty}\) is the asymptotic fraction of
  starch that can be hydrolysed by amylase; \(k_{0}\) is the initial
  reaction‑rate coefficient; \(k_{e}\) is the ``inherent'' reaction‑rate
  constant that reflects the intrinsic structure of the starch; \(\tau\)
  is a characteristic time constant that controls how quickly
  \(k(t)\) relaxes from \(k_{0}\) to \(k_{e}\). After roughly \(5\tau\)
  the coefficient has reached its asymptotic value.
\end{itemize}

\eqref{Eq.10b} explicitly describes why and how the reaction‑rate coefficient
changes over time. Physically, the model assumes that the rate
coefficient starts at the set‑up value \(k_{0}\) and then decays
exponentially toward the asymptotic value \(k_{e}\) with a first‑order
dynamics governed by \(\tau\).

\subsubsection{Analytical solutions}
\setcounter{equation}{10}
Solving the coupled system yields:
\begin{align}
    \begin{cases}
        X_{t} & = X_{\infty}\left( 1 - e^{- k_{e}t - \tau\left( k_{0} - k(t) \right)} \right), \\
        k_{t} & = k_{e} + \left( k_{0} - k_{e} \right)e^{- \frac{t}{\tau}}.
    \end{cases}
    \label{Eq.11}
\end{align}

If the initial and inherent rates are identical (\(k_{0} = k_{e}\)),
\eqref{Eq.11} reduces to the classical first‑order exponential model.

In practice, \(\tau\) can be treated as an independent parameter, but it
is expected to decrease as \(k_{0}\) increases, because a larger initial
rate accelerates the decay of \(k(t)\) toward \(k_{e}\). A useful
simplification is to set \(\tau = 1/k_{0}\), which eliminates one free
parameter and leads to a three‑parameter formulation:
\begin{align}
    \begin{cases}
X_{t} & = X_{\infty}\left( 1 - e^{- k_{e}t - \left( 1 - \frac{k(t)}{k_{0}} \right)} \right), \\
k_{t} & = k_{e} + \left( k_{0} - k_{e} \right)e^{- k_{0}t}.
\end{cases}
    \label{Eq.12}
\end{align}

\subsubsection{Calibration strategy}

Fitting candidate amylolysis models using only the integrated form
(i.e., the hydrolysis extent \(X_{t}\)) often yields satisfactory
goodness‑of‑fit scores for all models reported in the literature, but
this alone does not discriminate between them. A more rewarding 
strategy is to fit concomitantly the hydrolysis‑extent curve \(X_{t}\) and the
instantaneous reaction‑rate coefficient \(k_{a}(t)\) computed from
\eqref{Eq.3b)}, as described in the ``Extended procedure for
fitting of starch degradation models'' section. In this work, the four
models: first‑order \eqref{Eq.1}, Weibull \eqref{Eq.4}, the coupled system
\eqref{Eq.11} and its simplified version \eqref{Eq.12}, are calibrated against
both \(X_{t}\) and \(k_{a}(t)\) time-course measurements.

\begin{figure}
\centering
\includegraphics[width=\linewidth]{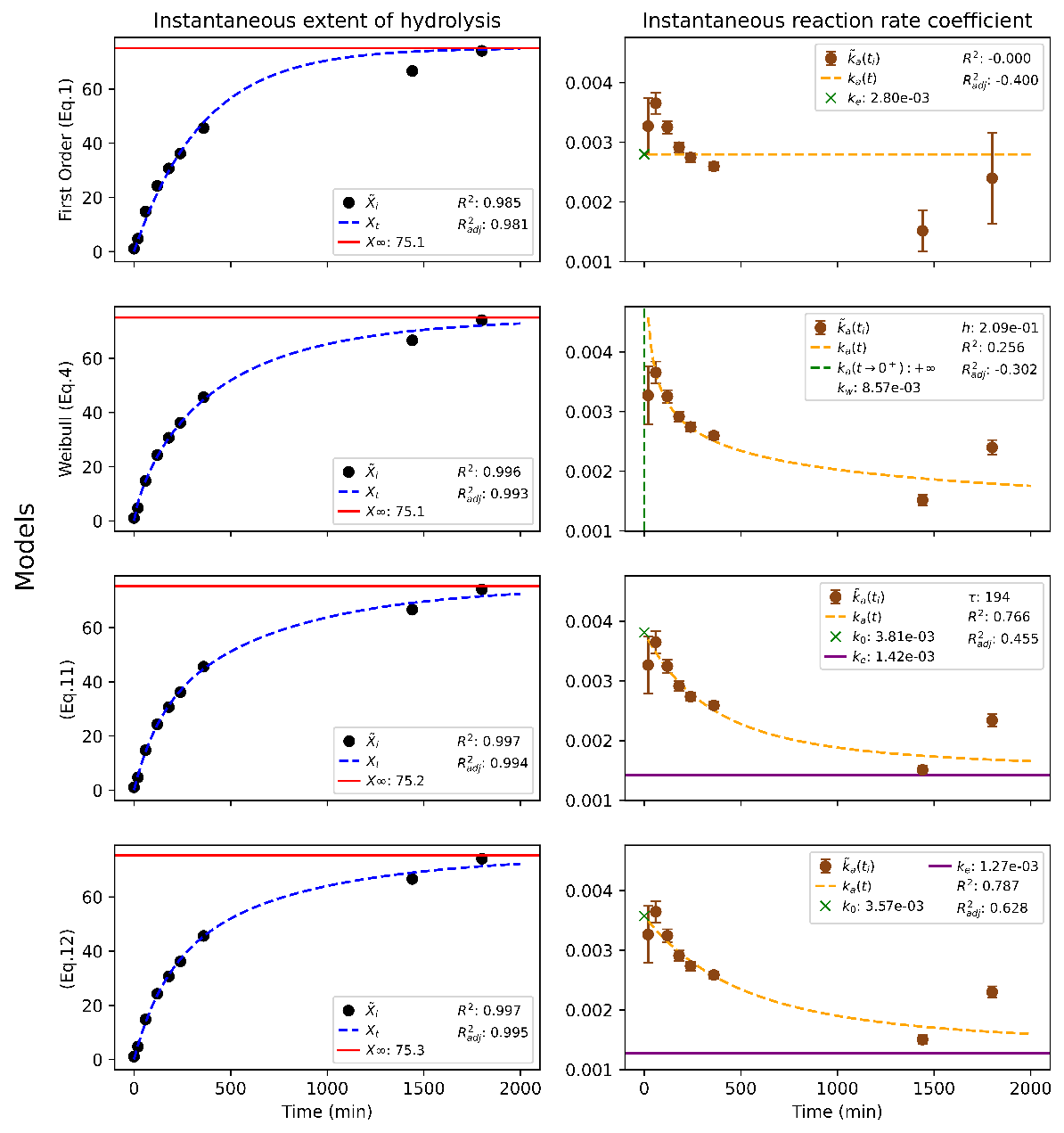}
\caption{Example of model fits to data from a WS sample. Each row
corresponds to a different model. The left panels show the calibration
based on the hydrolysis data (\(X_{t}\)), whereas the right panels show
the calibration based on the associated instantaneous reaction rate
($k_a(t)$).}
\label{Fig.5}
\end{figure}

Fig.\ref{Fig.5} shows the calibration of four candidate kinetic models for the
amylolysis of a single normal‑wheat‑starch (WS) sample. The fitting
procedure minimises the residuals for both the hydrolysis
extent\,\(X_{t}\) and the instantaneous reaction‑rate
coefficient\,\(k_{a}(t)\). It must be noted that \(k_{a}(t)\) 
values changes with the selected model along with the estimate of \(X_{\infty}\)
following \eqref{Eq.3c}.   
The quality of each fit is reported with two
statistics that appear on every panel: the coefficient of
determination\,\(R^{2}\) and the adjusted coefficient of
determination\,\(R_{\text{adj}}^{2}\), which penalises models that
contain more parameters (both calculated separately for\,\(X_{t}\)
and\,\(k_{a}(t)\)).

All four models reproduce the overall hydrolysis curve very well, as
indicated by\,\(R^{2}\) and\,\(R_{\text{adj}}^{2}\) values
exceeding $0.98$. This high level of agreement means that the extent data
alone cannot discriminate between the models. By contrast, the
right‑hand panels (which display the fit of the reaction‑rate data)
reveal clear differences. The two new models (Eqs.\,11 and\,12) achieve
substantially higher\,\(R^{2}\) values (0.766 and 0.787) than the
Weibull model (0.256) and the classic first‑order model (0.000). The
Weibull model's poor performance stems from its asymptotic behaviour at
early times
(\({\lim}_{t \rightarrow 0}k_{a}^{\text{Weibull}}(t) = \infty\)), a
feature not observed in the experimental data and therefore heavily
penalised in the\,\(R^{2}\) calculation.

Imposing the constraint \(\tau = 1/k_{0}\) (as in \eqref{Eq.12}) reduces the
number of free parameters and yields a higher\,\(R_{\text{adj}}^{2}\)
than \eqref{Eq.11}, reflecting the benefit of a more parsimonious
description. Appendix.D provides two additional examples of
fitting the digestion of a waxy‑wheat‑starch (WWS) 
and of a high‑amylose wheat starch (HAWS). The selected WWS kinetics illustrates 
the case of a \(k_{a}(t)\) growing over time. 

\begin{table}[]
    \centering
    \rowcolors{2}{black!10}{white}
    \begin{tabular}{c|cccc}
    \hline
    \multicolumn{5}{l}{1a) WS samples}\\
      &  First Order \eqref{Eq.1} & Weibull \eqref{Eq.4} & \eqref{Eq.11} & \eqref{Eq.12}\\
      \hline
    \cellcolor{white} \(X_{\infty}\) & \(75.07 \pm 3.96\) & \(75.25 \pm 4.10\) &
\(74.56 \pm 4.06\) & \(74.55 \pm 4.08\) \\
\(k_{e}\) \((10^{- 3})\) & \(2.60 \pm 0.29\) &
\(2.19 \in \lbrack 2.07,2.32\rbrack\)* & \(2.17 \pm 0.66\) &
\(2.27 \pm 0.78\) \\
\cellcolor{white} \(k_{0}\) \((10^{- 3})\) & & & \(2.74 \pm 1.08\) & \(2.82 \pm 0.45\) \\
\(k_{w}\) \((10^{- 3})\) & & \(3.91 \pm 1.61\) & & \\
\cellcolor{white} \(h\) & & \(0.062 \pm 0.073\) & & \\
\(\tau\) & {\color{white} First Order \eqref{Eq.1}} & {\color{white}Weibull \eqref{Eq.4}} & \(499 \pm 1288\) & \(364 \pm 58\) \\
\cellcolor{white} \(R^{2}\ (X)\) & \(0.996 \pm 0.004\) & \(0.997 \pm 0.002\) &
\(0.998 \pm 0.002\) & \(0.998 \pm 0.002\) \\
\(R_{adj}^{2}\ (X)\) & \(0.994 \pm 0.005\) & \(0.996 \pm 0.004\) &
\(0.996 \pm 0.004\) & \(0.997 \pm 0.003\) \\
\cellcolor{white} \(R^{2}\ (k_{a})\) & \(- 0.14 \pm 0.15\) & \(- 0.07 \pm 0.64\) &
\(0.21 \pm 0.46\) & \(0.14 \pm 0.46\) \\
\(R_{adj}^{2}\ (k_{a})\) & \(- 0.59 \pm 0.21\) & \(- 0.89 \pm 1.15\) &
\(- 0.89 \pm 1.13\) & \(- 0.54 \pm 0.88\) \\
\hline\hline
    \end{tabular}
    \rowcolors{2}{white}{black!10}
    \begin{tabular}{c|cccc}
\multicolumn{5}{l}{1b) WWS samples}\\
   \cellcolor{white}  \(X_{\infty}\) & \(85.60 \pm 4.13\) & \(85.56 \pm 4.16\) &
\(84.57 \pm 4.17\) & \(84.46 \pm 4.16\) \\
\(k_{e}\) \((10^{- 3})\) & \(5.14 \pm 0.84\) &
\(3.85 \in \lbrack 0.67,13.5\rbrack\)* & \(5.71 \pm 1.73\) &
\(6.65 \pm 1.35\) \\
\cellcolor{white} \(k_{0}\) \((10^{- 3})\) & & & \(3.02 \pm 2.36\) & \(4.76 \pm 1.24\) \\
\(k_{w}\) \((10^{- 3})\) & & \(4.32 \pm 2.25\) & & \\
\cellcolor{white} \(h\) & & \(- 0.056 \pm 0.084\) & {\color{black!10}\eqref{Eq.11}} & {\color{black!10}\eqref{Eq.12}} \\
\(\tau\) &{\color{white} First Order \eqref{Eq.1}} & {\color{white}Weibull \eqref{Eq.4}}& \(249 \pm 377\) & \(224 \pm 60\) \\
\cellcolor{white} \(R^{2}\ (X)\) & \(0.996 \pm 0.004\) & \(0.997 \pm 0.002\) &
\(0.998 \pm 0.002\) & \(0.994 \pm 0.004\) \\
\(R_{adj}^{2}\ (X)\) & \(0.987 \pm 0.008\) & \(0.989 \pm 0.007\) &
\(0.991 \pm 0.007\) & \(0.991 \pm 0.007\) \\
\cellcolor{white} \(R^{2}\ (k_{a})\) & \(- 0.43 \pm 0.14\) & \(- 0.84 \pm 0.69\) &
\(0.06 \pm 0.21\) & \(- 0.24 \pm 0.16\) \\
\(R_{adj}^{2}\ (k_{a})\) & \(- 1.00 \pm 0.19\) & \(- 2.22 \pm 1.21\) &
\(- 1.18 \pm 0.48\) & \(- 1.20 \pm 0.31\) \\
\hline\hline
\end{tabular}
    \rowcolors{2}{white}{black!10}
    \begin{tabular}{c|cccc}
\multicolumn{5}{l}{1c) HAWS samples}\\
   \cellcolor{white}  \(X_{\infty}\) & \(58.49 \pm 2.85\) & \(58.52 \pm 2.86\) &
\(58.21 \pm 2.82\) & \(58.28 \pm 2.82\) \\
\(k_{e}\) \((10^{- 3})\) & \(3.68 \pm 0.34\) &
\(2.04 \in \lbrack 0.61,5.46\rbrack\)* & \(1.71 \pm 0.42\) &
\(1.66 \pm 0.39\) \\
\cellcolor{white} \(k_{0}\) \((10^{- 3})\) & & & \(5.32 \pm 0.96\) & \(5.09 \pm 0.80\) \\
\(k_{w}\) \((10^{- 3})\) & & \(11.95 \pm 4.06\) & & \\
\cellcolor{white} \(h\) & & \(0.215 \pm 0.060\) & {\color{black!10}\eqref{Eq.11}} & {\color{black!10}\eqref{Eq.12}} \\
\(\tau\) &{\color{white} First Order \eqref{Eq.1}} & {\color{white}Weibull \eqref{Eq.4}}& \(200 \pm 109\) & \(201 \pm 32\) \\
\cellcolor{white} \(R^{2}\ (X)\) & \(0.984 \pm 0.009\) & \(0.995 \pm 0.003\) &
\(0.997 \pm 0.002\) & \(0.997 \pm 0.002\) \\
\(R_{adj}^{2}\ (X)\) & \(0.978 \pm 0.012\) & \(0.992 \pm 0.005\) &
\(0.994 \pm 0.004\) & \(0.995 \pm 0.003\) \\
\cellcolor{white} \(R^{2}\ (k_{a})\) & \(- 0.05 \pm 0.09\) & \(0.68 \pm 0.36\) &
\(0.85 \pm 0.26\) & \(0.80 \pm 0.31\) \\
\(R_{adj}^{2}\ (k_{a})\) & \(- 0.47 \pm 0.13\) & \(0.45 \pm 0.64\) &
\(0.65 \pm 0.61\) & \(0.66 \pm 0.54\) \\
\hline
    \end{tabular}
    \caption{Results of the digestion of WS (Table 1a), WWS (Table
1b), and HAWS (Table 1c) samples according to the different models,
based on the fitting of both \(X(t)\) and \(k_{a}(t)\). In the tables,
each quantity is summed up by its average ± its standard deviation. *:
estimations of \(k_{e}\) using the Weibull method, followed by their
95\% confidence intervals.}
    \label{tab:1}
\end{table}

\subsection{Uncertainty propagation for the instantaneous reaction‑rate
measurements}

In Fig.\ref{Fig.5} standard deviation are displayed on \(k_{a}(t)\) curves and
represent the uncertainty intervals that arise from propagating the
experimental error on the hydrolysis extent\,\({\widetilde{X}}_{i}\)
through \eqref{Eq.3b)}. If the standard deviation of the measured extent is
denoted \(\Delta X\), then the corresponding
uncertainty on the instantaneous rate at time\,\(t_{i}\) is

\[\Delta k_{a}(t_{i}) = \frac{\Delta X}{t_{i}(X_{\infty} - {\widetilde{X}}_{i})}.
\]

Assuming \(\Delta X\)
 is independent of time, it can be estimated from hydrolysis measurements 
at $t=0$. Indeed, since $\tilde{X}_0$, is supposed to be 0, then the standard error of 
hydrolysis measurements at $t=0$ is 
$\sqrt{\mathbb{V}({\widetilde{X}}_{0}^{2})}=\sqrt{\mathbb{E}({\widetilde{X}}_{0}^{2})}$,
which gives \(\Delta X = 0.67\) in this work.

This propagation is generally well‑behaved, but two limiting situations
merit special attention:

\begin{enumerate}
\def\labelenumi{\arabic{enumi}.}
\item
  \textbf{Early‑time points (\(t \rightarrow 0\)\textbf{)}:} when the
  denominator contains a very small\,\(t\) values, \eqref{Eq.3b)} becomes highly
  sensitive to measurement noise, so the derived\,\(k_{a}(t)\) values
  are unreliable.
\item
  \textbf{Late‑time points
  (\({\widetilde{X}}_{i} \rightarrow X_{\infty}\)\textbf{)}:} if
  \({\widetilde{X}}_{i} \geq X_{\infty}\)the expression
  for\,\(k_{a}(t)\) is undefined; when
  \({\widetilde{X}}_{i} = X_{\infty} - \varepsilon\) with
  \(\varepsilon \ll 1\), the uncertainty blows up as
  \(\Delta k_{a}(t_{i}) = \frac{\Delta X}{t_{i}\varepsilon} \underset{\epsilon\to 0}{\rightarrow} \infty\)
  even for "long" $t_i$ encountered during experiments as $\widetilde{X}_{i}$ tends to $X_\infty$
  exponentially in time. Consequently, when hydrolysis has reached its asymptotic rate, 
  \(k_{a}(t)\) estimates carry very large confidence intervals. 
\end{enumerate}

High uncertainty for \(k_{a}(t)\) for late-time points is reflected by 
large variations of \(k_{a}(t)\) values when \(\widetilde{X}_{i}\) gets close to \(X_{\infty}\) 
Paradoxically, models with large overestimates of \(X_{\infty}\) are therefore, less affected and might obtain
better fitting scores.
Obtaining an independent measurement of \(k_{a}(t_{i})\) would reduce the propagation of
uncertainty from \eqref{Eq.3b)} and give a more reliable fitting results.

\subsection{Summary of the model‑fitting results (Table\ref{tab:1})}

For the 652 normal wheat‑starch (WS), 14 waxy wheat‑starch (WWS) and 14
high‑amylose wheat‑starch (HAWS) samples, all four kinetic models
provide an excellent description of the hydrolysis extent\,\(X_{t}\)
(average\,\(R^{2} > 0.99\)). The adjusted\,\(R^{2}\) values are very
close to each other, so the digestogram alone does not allow a clear
ranking of the models, except for the first‑order model on the HAWS
data, whose performance is slightly lower.

The fitted asymptotic hydrolysis fractions\,\(X_{\infty}\) are
consistent across models: $\approx$\,75\,\% for WS, $\approx$\,85\,\% for WWS and
$\approx$\,58\,\% for HAWS.

The ``inherent'' rate constant, \(k_{e}\) obtained with the two new
models (Eqs.\,11\,and\,12) is similar for WS ($\approx 2.2\times 10^{-3}$\,min$^{-1}$)
and HAWS ($\approx 1.7\times 10^{-3}$min$^{-1}$). The classic first‑order model yields
 higher average values for \(k_{e}\) ($\approx 2.6\times 10^{-3}$min$^{-1}$ for WS and
$\approx 3.7\times 10^{-3}$min$^{-1}$ for HAWS). For the Weibull model,\,\(k_{e}\) is
not a direct parameter; it can be estimated a posteriori using the
linear relationship between\,\(k_{w}\) and\,\(h\) (see section \ref{revealing ke}). 
This estimate is close to fitted \(k_{e}\) from
(Eqs.\,11\,and\,12) for WS, whereas for WWS and HAWS the confidence
intervals are too broad to draw firm conclusions.

The initial rate\,\(k_{0}\) (Eqs.\,11\,and\,12) is on average slightly higher than the
fitted\,\(k_{e}\) for WS, much higher for HAWS, and lower for WWS. 
The results for HAWS and WS are consistent with the data (Fig.\ref{Fig.1}b,d) despite
large standard deviations. The results for WWS, on the other hand, are inconsistent
with Fig.\ref{Fig.1}c, which shows that the two models are inadequate to describe this 
type of kinetics. The characteristic time\,\(\tau\) is similar in both models; however,
\eqref{Eq.12} (where \(\tau\) is fixed to\,\(1/k_{0}\)) shows a much reduced
standard deviation compared to \eqref{Eq.11}, reflecting the reduced degrees of freedom.
The large standard‑deviations obtained \(\tau\) fitted with \eqref{Eq.11}, do not allow
to comment on the duration of the transition from\,\(k_{0}\) to\,\(k_{e}\). In \eqref{Eq.11}
model, \(\tau\) is certainly a parameter not correctly fitted and possibly useless in a number
of cases.

With respect to the Weibull parameters\, fitted values of \(k_{w}\) for WS and HAWS 
appear to be on average higher than \(k_{0}\). In theory \(k_{w}= k(t=1)/(1-h)\) 
\citet{Kansou2015a}, therefore \(k_{w}\) should be equal to \(k_{0}\) when \(h\) 
is null and lower than \(k_{0}\) when \(h\) is positive. This indicates an overestimation 
which stems from the asymptotic behaviour of the Weibull when time is close to zero. 
The parameter\,\(h\) is on average positive for HAWS (exponential‑type saturation), 
and close to zero for WS and WWS (with high standard-deviation), which seems consistent
with the data (Fig.\ref{Fig.1}). Note that the values of \(k_{w}\) and \(h\) plotted in 
Fig.\ref{Fig.3}, were obtained through the fitting of only \(X_{t}\)
for the sake of the demonstration, while the values reported in Table\ref{tab:1}
results from dual fitting of \(X_{t}\) and \(k_{a}(t)\). This explains the observed 
discrepancies. 

Overall, the combined analysis of hydrolysis extents and instantaneous
rates, together with the uncertainty assessment, provides a robust basis
for discriminating kinetic models and for estimating the intrinsic
reaction‑rate constant\,\(k_{e}\).

\begin{figure}
\centering
\includegraphics[width=\linewidth]{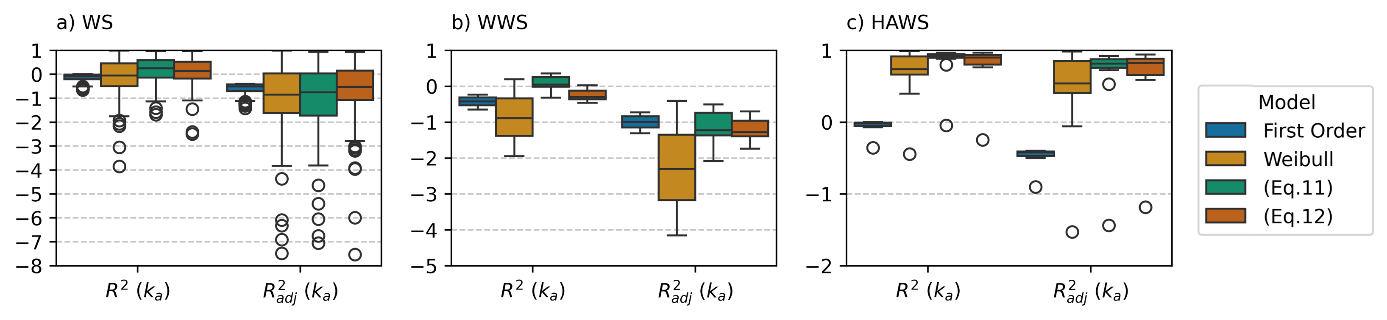}
\caption{Boxplots of \(R^{2}\) and \(R_{adj}^{2}\) scores' distribution
associated with the instantaneous reaction rate coefficient (\(k_{a}\))
fits for WS (6a), WWS (6b), and HAWS (6c) samples, distinguished by
models: First Order \eqref{Eq.1} in blue, Weibull \eqref{Eq.4} in orange, \eqref{Eq.11} in
green, and \eqref{Eq.12} in red. Each box shows the quartiles of the dataset,
with a horizontal bar at the median level, while the whiskers extend to
show the rest of the distribution, except for points that are samples
considered as outliers.}
\label{Fig.6}
\end{figure}

The goodness‑of‑fit for the instantaneous reaction‑rate
coefficient\,\(k_{a}(t)\) varies markedly across the four kinetic models
and the three starch types (WS, WWS, HAWS). In addition to the mean
and standard‑deviation values of the \(R^{2}\) and adjusted\,\(R^{2}\)
scores reported in Table\ref{tab:1}, we visualised the full distribution of
these scores in Fig.\ref{Fig.6} to enable a more detailed
discrimination of the models. Some readers may be surprised
by the low values of \(R^{2}\) obtained for \(k_{a}(t)\).
It is well known that this criteria is strongly affected 
by the presence outliers. In this work they arise from the 
application of \eqref{Eq.3b)} to late-time points.
An illustrative example can be found in Appendix.D (Fig.9) for the fitting of \(k_{a}(t)\)
with \eqref{Eq.11}. While the corresponding \(R^{2}\) is 0.24, removing the outlier observed 
at t=1440 min is sufficient to get a \(R^{2}\) of nearly 0.99. For this reason absolute \(R^{2}\)
values should be considered with care.

\begin{itemize}
\item
  \textbf{First‑order model \eqref{Eq.1}:} this model shows the smallest
  variability in scores, but the \(R^{2}\) values are close to zero for
  WS and HAWS ($\approx$\,0) and negative for WWS ($\approx$\,‑0.4), indicating a poor
  description of the \(k_{a}(t)\) dynamics.
\item
  \textbf{Weibull model \eqref{Eq.4}:} the \(R^{2}\) scores are also near
  zero for WS, strongly negative for WWS ($\approx$\,‑0.8), and moderately
  positive for HAWS ($\approx$\,0.68). Thus, the Weibull model captures the
  dynamics only for the high‑amylose samples.
\item
  \textbf{Model based on \eqref{Eq.11}:} on average this model yields the
  highest \(R^{2}\) values, yet the scores remain low for WS ($\approx$\,0.21)
  and WWS ($\approx$\,0.06) and are only satisfactory for HAWS, where only a
  single sample shows a poor fit.
\item
  \textbf{Model based on \eqref{Eq.12}:} the performance is comparable to
  \eqref{Eq.11} but the average \(R^{2}\) is slightly lower; the pattern of
  results across the three starch types is the same.
\end{itemize}

When the adjusted\,\(R^{2}\) scores are examined (Fig.\ref{Fig.6}, right‑hand
panels), models' performances depend on the type of wheat samples:

\begin{itemize}
\item
  \textbf{WS:} over the large diversity of time-course measurements 
  displayed by the hydrolyses of 652 WS, the two‑parameter first‑order model provided
  good fitting in general, but the other models performed better in about 
  two-thirds of cases. Unfortunately they also obtained very negatives \(R^{2}\) scores 
  in several occasions which impaired their overall fitting performance, while the first‑order model
  proved more robust. Overall the three‑parameter model of \eqref{Eq.12} achieved the 
  highest average adjusted\,\(R^{2}\), reflecting a better balance between fit quality 
  and model complexity.
  
\item
  \textbf{WWS:} the first‑order model obtains the
  highest average adjusted\,\(R^{2}\) because models with more parameters 
  failed to reproduce the non‑monotonic behaviour of \(k_{a}(t)\);
  consequently, the penalty for extra parameters lowers their adjusted
  scores. The Weibull model consistently yielded the lowest
  adjusted\,\(R^{2}\) for this starch type.
\item
  \textbf{HAWS:} the first‑order model is overly simplistic and ranks
  as the worst kinetic model. The Weibull model is a viable alternative but
  showed a broad adjusted‑score distribution. Both \eqref{Eq.11} and
  \eqref{Eq.12} fit the \(k_{a}(t)\) kinetics very well, with \eqref{Eq.11}
  displaying a narrower distribution than \eqref{Eq.12}.
\end{itemize}

In Summary, the two new formulations \eqref{Eq.11} and \eqref{Eq.12}
outperform the first-order model and the Weibull model in most cases. 
The choice between them depends on the starch type: \eqref{Eq.11} is 
preferable for HAWS, whereas \eqref{Eq.12} gives the best overall 
performance for WS. However, none of the current models captures 
the non‑monotonic \(k_{a}(t)\) patterns observed for WWS, indicating 
that the differential equation governing the rate coefficient 
\eqref{Eq.10b} should be refined to accommodate such behaviour. 

\section{Discussion}

Time-course measurements of hydrolysis extent of starch are commonly fitted with kinetics 
equations as a part of the data analysis process. The choice of the math model
depends on the purpose of the study. For example analysing the action of amylase is usually 
performed by means of the Michaelis-Menten kinetics model, whose parameters have clear 
physical meanings. Unfortunately, the Michaelis-Menten equation does not generally fit 
well the starch hydrolysis curve beyond a short period \citet{Butterworth2022}. 
Thus, fitting the complete digestion of starch \emph{in vitro} is generally performed by means 
of empirical models. Various empirical models have been proposed to investigate the starch-amylase system 
over time of hydrolysis \citet{Nguyen2018}. Although they are empirical, every meaningful 
models assume a theoretical physical structure and a mode of operation of the system under study. 
The fitting of experimental data indicates whether or not the underlying assumption is plausible.
Thus it is well established in the domain that starch digestion kinetics 
apparently matches a first-order reaction kinetic. So not surprisingly, popular empirical model
are variants of the first-order reaction kinetic model.

In this work we have shown that different empirical models can all obtain high goodness-of-fit for 
starch hydrolysis extent, $X_t$. Therefore, fitting the hydrolysis extent alone proved to be not 
sufficient to discriminate the models and to gain a better understanding of the starch-enzyme 
interaction over time. Other criteria should then be examined to select a model. 
\citet{Butterworth2012} proposed to use a Log of Slope (LOS) plot, i.e. $ln(dX/dt)$ 
against $t$, mainly to improve the determination of $X_\infty$ but also to track changes 
in $k_a(t)$. Indeed, LOS plot is linear of slope $k_a(t)=-k$ when first-order kinetics 
applies. Non-linear LOS plot curves are generally interpreted as two starch fractions digested 
in a row at different rates, which leads to a sequential two-exponential model 
\citet{Edwards2014}. Plotting the instantaneous rate coefficient following \eqref{Eq.3b)}, 
is comparable to using LOS plot to assess the $k_a(t)$. In agreement with the present work, 
the majority of non-linear LOS plots of starch digestion reported in the literature display 
a decrease of $k_a(t)$ with time that converges toward a constant value (e.g. in \citet{Edwards2014})
which would typically correspond to the situation of HAWS hydrolysis (Fig.\ref{Fig.1}d). The present 
work goes a step further by combining the fitting of the hydrolysis extent $X_t$ and the instantaneous 
rate coefficient $k_a(t)$, using an overall fitting criterion.
Thus, a direct strategy was chosen to obtain estimates of $k_a(t)$ 
for any time-point which involves the concomitant estimation of $X_\infty$. This strategy 
is more selective and hopefully meaningful for in-depth starch digestion kinetics analysis. 
Unfortunately it produces highly uncertain estimates when $X_t$ gets close to $X_\infty$, as reflected 
by very low \(R^{2}\) scores. In this study, the large number of wheat starch hydrolysis runs
made it possible to obtain consistent results nonetheless. The LOS plot could be used 
to assess $k_a(t)$ independently of $X_\infty$, even though it has been also recommended to exclude
the experimental points from the final stage of hydrolysis in LOS plot analysis \citet{Butterworth2012}.

In addition to assess the fitting performance of the candidate models, it is also important to compare 
their characteristics and their properties from a mathematical and physical perspective. 
We suggest the following set of criteria : minimum number of parameters, continuous and 
differentiability of class $\mathcal{C}^{\mathbf{1}}$ or better, and for physical interpretation: being equivalent 
to the classical exponential Eq.1 for a subset of parameters, being flexible enough to fit 
rapid-to-slow pattern and slow-to-rapid pattern, present finite limits for $k_a(t)$ when 
$t\to0$ and $t\to\infty$. 

The model flexibility, which in the context of starch digestion kinetics represents the 
capacity of the model to fit saturation-like curve or S-shape curve, seems particularly 
useful for nutritional studies, where \emph{in vitro} digestion assays are routinely employed to 
estimate \emph{in vivo} glycaemic indices. For large food particles or starch embedded in plant 
tissues, the diffusion of the enzyme is expected to be slower and the digestion rate lower, 
which generates a kinetic profile different than for purified starch digestion \citet{Butterworth2022}. 
For example, \citet{PallaresPallares2018} reported pronounced S-shaped digestibility 
curves for starch in bean cotyledon cells, with lag phases up to 60 min caused by cell wall 
barriers that impede enzyme diffusion. Those authors fitted the data with a logistic model, 
which is of course fine for a descriptive analysis. However, the models explored in this 
paper, for example the Weibull model or \eqref{Eq.11} and \eqref{Eq.12} models, could have 
provided more insights about the increase of $k_a(t)$ over time and then about the status 
of the starch-enzyme system. In other terms, they would make it possible to quantify the 
impact of the surrounding matrix on the initial deviation from first-order reaction kinetics. 

At the light of this work, presenting finite and meaningful limits for $k_a(t)$ appears to 
be desirable criterion for a kinetics model of starch digestion. Analysis of hundreds of time 
course measurements on wheat starches hydrolysed by amylase showed that the instantaneous 
reaction rate coefficient $k_a(t)$ tends to converge to a constant value, $k_e$. 
A formal proof of the existence of this constant is given in the supplementary material 
using the Weibull model (see Appendix.B). It appears that $k_e$ reflects an intrinsic 
property of the starch–enzyme system and represents a stable equilibrium from a 
dynamical systems perspective.
The estimated $k_e$ is likely to be specific to each starch type and then mirrors their 
inherent susceptibility to amylase, even though an in-depth analysis would be required 
to strengthen this claim. By contrast, the initial value of $k_a(t)$, denoted $k_0$, is 
highly variable and must primarily depend on surface-related factors which in turn are 
more sensitive to the environment of the experiment.
An average estimates of $k_0$ can be obtained from several replicates at the condition 
to fit each replicate separately. This is the approach we used in this work to discriminate 
between properties sensitive to external factors from inherent and stable properties of the system. 
Both parameters, $k_0$ and $k_e$, are however important for describing starch digestion 
kinetics and should be explicit in the kinetic model used to fit experimental data 
whenever it is possible.
A comparison of kinetics models following the criteria discussed above is supplied in Table \ref{tab.2},
from Appendix.A. 
It includes the four models investigated in this work and two types of double exponential 
models, parallel and sequential, that are commonly used in starch digestion studies to fit 
time-course measurements. This table reveals the limitations of the various models. The 
traditional first-order model remains the most widely used in the literature for good 
reasons. Our results show that with the minimum number of parameters it captures the 
overall behavior of the system reasonably well, although it is overly simplistic when $k_a(t)$ 
changes with time. Double exponential models assume distinct starch fractions each following its 
own first-order kinetics and distinct susceptibility levels. The parallel model would postulate 
two different types of starches \citet{Li2019}, which is generally not justified when only 
one type starch is being hydrolyzed. The physical interpretation of the sequential model reflects 
a mechanism that distinguishes between the hydrolysis of accessible starch at the substrate surface 
and inside starch \citet{Edwards2014}, which seems more consistent with the observations. 
However the discontinuity it introduces at an arbitrary time-point does not seem justified. 
From a mathematical point of view the parallel model appears more relevant than the sequential 
which is not of class $\mathcal{C}^{\mathbf{1}}$ at $t=T$. Double-exponential models have at least 
four parameters and they are not flexible, which are strong limitations. Additionally, we noticed 
numerically that both double exponential models exhibit at least one degree of freedom within parameters. 
Indeed, despite very good fitting scores for $X_t$ and $k_a(t)$ curves, $X_\infty$ is often naturally 
calibrated to 100\% (Fig.\ref{Fig.7}). Numerical constraints could be added to correct this behavior, but this
shows that those models are over-parameterized for the starch digestion kinetics.
The Weibull model main drawback is the physical interpretation, as it does not have finite 
limits for $k_a(t)$ in the general case. In practice the interpretation of $k_w$ and $h$ is 
not straightforward. To our knowledge, none of the common empirical models used to fit 
in vitro starch digestion kinetics satisfies all the criteria we listed. To fill this gap, 
we introduced a four parameter model \eqref{Eq.11} and a three parameter model \eqref{Eq.12}. 
This type of model corresponds to a continuous density gradient in the hydrolysis rate between
the surface of the particles and their core, captured by the parameters $k_0$ and $k_e$. 
These models are consistent with a view where the early phase of hydrolysis is controlled 
by surface related factors of the substrate (specific surface area, degree of order of the 
outer glucan chains, the presence of surface proteins and lipids…), whereas the later phase 
of hydrolysis is governed by the organization and composition of the bulk starch matrix and 
can be described by a first order kinetics, as put forward in several publications 
(e.g. \citet{Warren2011,Butterworth2012,Wang2023}).
For the kinetics that match this pattern the models perform very well, for example for the 
digestion of HAWS (Fig.\ref{Fig.1}d). For non-monotonics kinetics of $k_a(t)$ as for WWS 
digestion (Fig.\ref{Fig.1}c), a more sophisticated model would be needed or the usual 
first-order model should be used otherwise. 
In practice a more reliable experimental procedure for the measurement of $k_a(t)$ would 
be needed to assess further the starch digestion kinetics model in the future. 

\section{Conclusion}

We have developed a new method for fitting and analysing in vitro starch digestion kinetics 
that combines hydrolysis extent over time with the instantaneous reaction-rate coefficient. 
The method exploits further the time-course measurements and should help to improve the kinetic
analysis of the hydrolysis of starch by amylase. Using a large dataset of starch digestion 
time-course measurements, we provided evidences that the reaction rate coefficient converges
over time towards a constant $k_e$, in which case the kinetics reduces to an apparent first-order
reaction kinetics. This constant might be a stable characteristic of starch-amylase 
system and reflect the starch macromolecular composition and organisation. At the early stage of
the hydrolysis, the system dynamics is unstable and exhibits deviations from the first-order reaction. 
The deviations can arise from the substrate surface conditions or from other uncontrolled factors. 
The variability observed at the early stages generally diminishes with starch degradation as the 
starch bulk matrix becomes dominant. 
We provide a procedure for estimating $k_e$ and introduce new kinetics models more consistent with these 
observations and with the current understanding of starch digestion mechanism. The findings and 
propositions provided in this paper should favor the analysis of the multiple factors 
impacting the starch digestion kinetics.









\begin{landscape}
\appendix
\section{Comparison of usual and new models' mathematical expression to fit
  hydrolysis data}\label{Appendix.A}
\begin{table}[H]
\centering
    \rowcolors{2}{white}{black!10}
    \begin{tabularx}{\linewidth}{>{\raggedright\arraybackslash}p{3cm}|*{6}{>{\centering\arraybackslash}X}}
    \hline
    Models & First order & Parallel \citet{Li2019} & Sequential \citet{Edwards2014} & Weibull & \eqref{Eq.11} & \eqref{Eq.12} \\
    \hline
    Parameters & \(X_{\infty}\), \(k_{e}\) & \(X_{\infty}\),\(\ k_{1}\),
\(k_{2}\), \(\alpha<X_\infty\) & \(X_{\infty}\), \(k_{1}\), \(k_{2}\), \(T\), $Y$ &
\(X_{\infty}\), \(k_{w}\), \(h < 1\) & \(X_{\infty}\), \(k_{e}\),
\(k_{0}\), \(\tau\) & \(X_{\infty}\), \(k_{e}\), \(k_{0}\) \\
\(X_t\) &
\(X_{\infty}(1 - e^{- k_{e}t})\) &
\(\alpha\left( 1 - e^{- k_{1}t} \right) + {(X}_{\infty} - \alpha)(1 - e^{- k_{2}t})\)
& \(Y\left( 1 - e^{- k_{1}t} \right),\) if \(t \leq T\),

\({X(T) + (X}_{\infty} - X(T))(1 - e^{- k_{2}(t - T)})\), if \(t > T\).
& \(X_{\infty}(1 - e^{- k_{w}t^{1 - h}})\) &
$X_{\infty}( 1 $-$e^{- k_{e}t - \tau\left( k_{0} - k(t) \right)})$
&
$X_{\infty}( 1 $-$e^{- k_{e}t - \left( 1 - \frac{k(t)}{k_{0}} \right)} )$ \\
$X$ continuous? class \(\mathcal{C}^{\mathbf{1}}\)? & yes, yes & yes, yes &
yes, {\color{red}no} & yes, yes & yes, yes & yes, yes \\
How to match the First order model? & - & \(k_{1} = k_{2}\) & \(k_{1} = k_{2}\), and $Y=X_\infty$ & \(h = 0\) &
\(k_{0} = k_{e}\) & \(k_{0} = k_{e}\) \\
\(k_a(t)\) & \(k_{e}\) &
\(- \frac{1}{t}ln(\frac{\alpha}{X_{\infty}}e^{- k_{1}t} + (1 - \frac{\alpha}{X_{\infty}})e^{- k_{2}t}))\)
& \(-\frac{1}{t}ln(1-\frac{Y}{X_\infty}(1-e^{-k_1 t}))\), if \(t \leq T,
k_{2}\frac{t-T}{t}-\frac{1}{t}ln(1-\frac{X(T)}{X\infty})\), if \(t > T\) & \(k_{w}(1 - h)t^{- h}\) &
\(k_{e} + \tau(k_{0} - k_{e})\frac{1 - e^{- \frac{t}{\tau}}}{t}\) &
\(k_{e} + (1 - \frac{k_{e}}{k_{0}})\frac{1 - e^{- k_{0}t}}{t}\) \\
\(\underset{\mathbf{t \rightarrow 0}}{\lim}{k_a(t)}\)
& \(k_{e}\) &
\(\frac{\alpha}{X_{\infty}}k_{1} + \left( 1 - \frac{\alpha}{X_{\infty}} \right)k_{2}\)
& \(\frac{Y}{X_\infty}k_{1}\) & \(\left\{ \begin{aligned}
{\color{red}0,\ \ } &{\color{red} h < 0} \\
k_{w},\ \  & h = 0 \\
{\color{red} + \infty,\ \  }& {\color{red}h > 0}
\end{aligned} \right.\ \) & \(k_{0}\) & \(k_{0}\) \\
\(\underset{\mathbf{t \rightarrow + \infty}}{\lim}{k_a(t)}\)
& \(k_{e}\) & \(min(k_{1},k_{2})\) & \(k_{2}\) &
\(\left\{ \begin{aligned}
{\color{red} + \infty,\ \  }& {\color{red}h < 0} \\
k_{w},\ \  & h = 0 \\
{\color{red}0,\ \  }& {\color{red}h > 0}
\end{aligned} \right.\ \) & \(k_{e}\) & \(k_{e}\) \\
Model rapid-to-slow pattern? ($X''_t<0$) & {\color{red}no} & yes & yes & yes, case: \(h > 0\) &
yes, case: \(k_{0} > k_{e}\) or $\frac{k_e-k_0}{\tau}< k_0^2$ & yes, case: \(2 k_{0} > k_{e}\) \\
Model slow-to-rapid pattern? (S-shape) & {\color{red}no} & {\color{red}no} & {\color{red}no} & yes, case: \(h < 0\) &
yes, case: \(k_{0} < k_{e}\) and $\frac{k_e-k_0}{\tau}>k_0^2$ & yes, case: \(2k_{0} < k_{e}\) \\
    \end{tabularx}
    \caption{Mathematical expressions and properties of usual (First order, Double exponential, Sequential, Weibull) and new 
    (\eqref{Eq.11}, \eqref{Eq.12}) models regarding the fit of hydrolysis data. Cells with red text show criteria that are not suitable in a model, regarding physical behaviors expected from literature. All parameters,
    except $h$, are supposed to be strictly non-negative.}
    \label{tab.2}
\end{table}
\end{landscape}

\begin{figure}
\centering
\includegraphics[width=\linewidth]{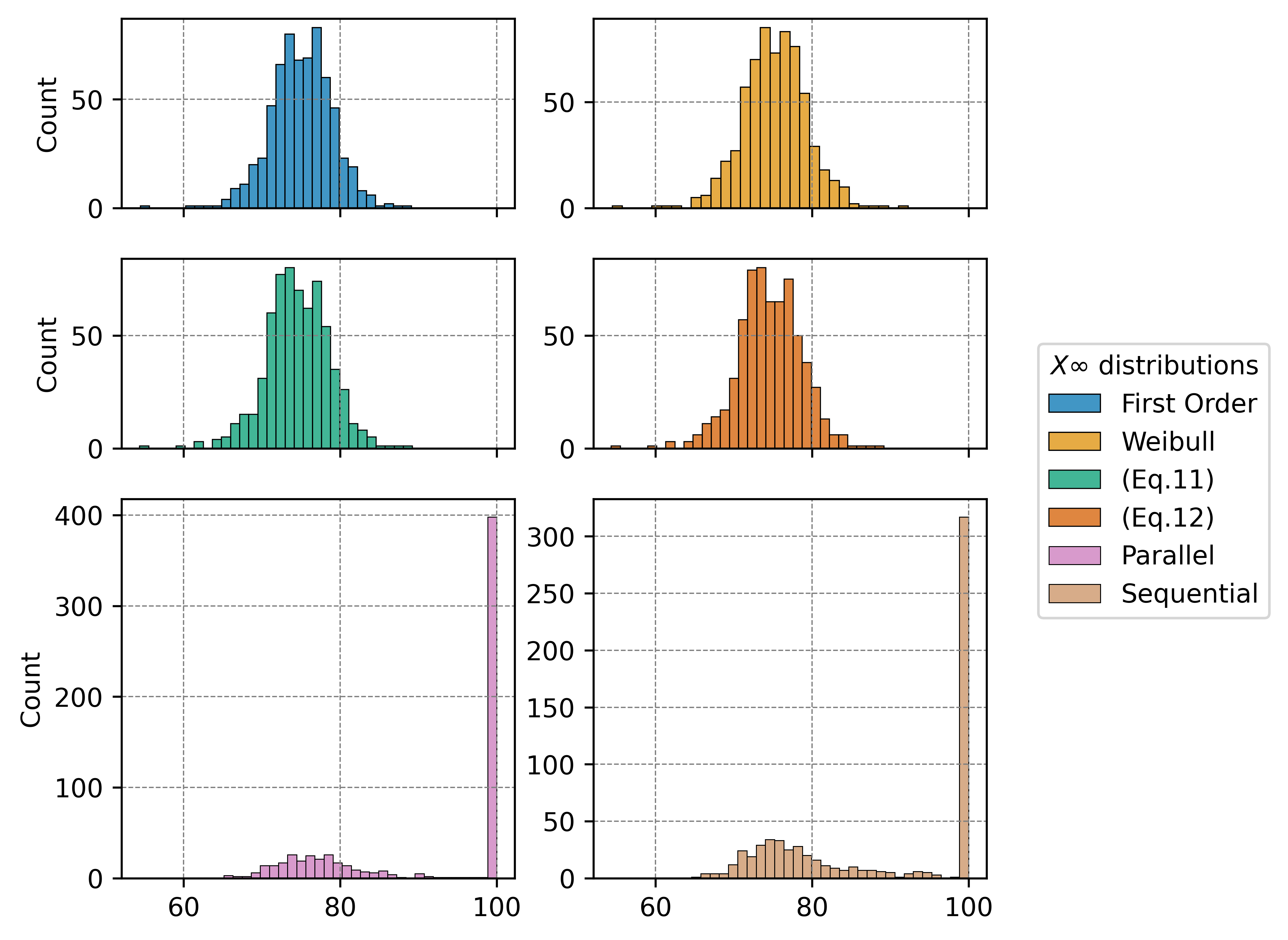}
\caption{Histograms of $X_\infty$ distributions obtained for WS samples depending on the models 
available in Table \ref{tab.2}.}
\label{Fig.11}
\end{figure}

The parallel \citet{Li2019} and sequential \citet{Edwards2014} models, 
introduced in Table \ref{tab.2}, are examples of hydrolysis models with double 
exponential expressions, i.e., considering that there are in fact two dynamics, 
in parallel or sequential, with constant reaction-rate coefficients $k_1$ and $k_2$,
rather than having a reaction-rate coefficient depending on time.
These models have a simple expression and mathematical interpretation for $X_t$,
even showing strong fitting capabilities on $X_t$ and $k_a(t)$.
However, as shown in Table\ref{tab.2}, those models are not physically validated, 
and compared to the one of $X_t$, the expressions of $k_a(t)$ are not easily interpretable
regarding the parameters of the models.
Additionally, in Fig.\ref{Fig.11}, we inspect the distribution of 
$X_\infty$ across all models from Table \ref{tab.2} on WS samples. Only the parallel 
and sequential models display $X_\infty$ fits at 100\%, even though they do fit both $X_t$
and $k_a(t)$ curves well. This is likely due to at least a degree of freedom in 
parameter choice which can compensate the fit of $X_\infty$ onto another parameter, since
the double exponential models have a significant number of parameters 
(4 and 5 respectively).

\appendix
\setcounter{section}{1}
\setcounter{equation}{0}
\renewcommand{\theequation}{B.\arabic{equation}}
\section{Proof of existence of \(k_{e}\)}
\label{Appendix.B}
Previous applications of the Weibull model \eqref{Eq.4} to fit \emph{in vitro}
digestions of maize and wheat starches showed that \(X_{\infty}\) was
not correlated with and \emph{k\textsubscript{w}} nor with \emph{h},
whereas \emph{k\textsubscript{w}} and \emph{h} were found to be highly
and positively correlated
\citet{Olawoye2020,Wang2022a}.
Wang et al., (2023), also shown that \emph{k\textsubscript{w}} and
\emph{h} were related to surface-related features, e.g. granule size,
while $X_\infty$ were clearly related to the starch
matrix composition and structure, e.g. chain length distribution, type-B
crystallinity and so on.

Fig.\ref{Fig.3} shows that even though the Weibull parameters \(k_{w}\) and \(h\)
are not theoretically linked together, there exists a logarithm
correlation between them when the Weibull model is fitted on hydrolysis
curves. This relationship is given by Equation \eqref{(Eq.B1)}, where \(a\) and
\(b\) are constants.

\begin{equation}
    h = a \times \ln\left( k_{w} \right) + b.
    \label{(Eq.B1)}
\end{equation}

Such relationship cannot be coincidental, and we demonstrate that it
proves the existence of a limit value of the instantaneous rate
coefficient, denoted \(k_{e}\), common to starch with identical origins.

Numerically, the fitting of the Weibull model is performed using a
classical least square approach with hydrolysis measurements, i.e.,
minimizing the error between the theoretical hydrolysis extent at time
\(t_{i}\), which is
\(X\left( t = t_{i} \right) = X_{\infty}(1 - e^{- k_{w}t_{i}^{1 - h}})\),
and its measured value \({\widetilde{X}}_{i}\). For simplicity of the
demonstration, let's assume that we calibrate \(k_{w}\) and \(h\) by
minimizing the sum of the squared residuals:

\begin{equation}
    J\left( k_{w},h \right) = \ \frac{1}{2}\sum_{i}^{}\left( k_{w}t_{i}^{1 - h} - {\widetilde{y}}_{i} \right)^{2},
    \label{(Eq.B2)}
\end{equation}

where
\(\widetilde{y_{i}} = \  - \ln\left( 1 - \frac{\ {\widetilde{X}}_{i}}{X_{\infty}} \right)\)
, which is also equal to
\({\widetilde{k}}_{ap}(t = t_{i}) \times t_{i}\) according to (Eq.3b).

Equation \eqref{(Eq.B2)} comes from rewriting
\(\frac{1}{2}\sum_{i}^{}\left( X(t = t_{i}) - {\widetilde{X}}_{i} \right)^{2}\),
by isolating the terms in the exponential expression. Doing so
simplifies analytical expressions of the partial derivatives of the
residuals, and is convenient to prove the existence of \(k_{e}\).

Minimization of \(J\) implies that
\(\frac{\partial J}{\partial k_{w}}\left( k_{w},h \right) = \frac{\partial J}{\partial h}\left( k_{w},h \right) = 0\),
which gives

\begin{align}
    \begin{cases}
        \frac{\partial J}{\partial k_{w}}\left( k_{w},h \right) & = \sum_{i}^{}{t_{i}^{1 - h}\left( k_{w}t_{i}^{1 - h} - {\widetilde{y}}_{i} \right)} = 0, \\
        \frac{\partial J}{\partial h}\left( k_{w},h \right) & = - \sum_{i}^{}{\ln\left( t_{i} \right)t_{i}^{1 - h}\left( k_{w}t_{i}^{1 - h} - {\widetilde{y}}_{i} \right)} = 0.
    \end{cases}
    \label{(Eq.B3)}
\end{align}

Both expressions in Equation \eqref{(Eq.B3)} are weighted sum of the residuals.
With the Weibull model, the weighting, \(t_{i}^{1 - h}\) in the first
expression and \(\ln\left( t_{i} \right)t_{i}^{1 - h}\) in the second,
is not uniform in time but increasing with \(t_{i}\) due to the modeling
requesting \(h < 1\). Thus, if there is not a perfect match of parameter
for the fitting, i.e., all residuals are not all close to zero, the
residual associated with the highest time measurement, denoted
\(t_{\max}\), will be the major term enforced to be close to 0. This
means that for each sample, numerically calibrating \(k_{w}\) and \(h\)
will enforce

\begin{equation}
    k_{w}t_{\max}^{1 - h} - {\widetilde{y}}_{\max} \approx 0\ \overset{\Leftrightarrow}{}k_{w}t_{\max}^{1 - h} - {\widetilde{k}}_{ap}(t = t_{\max})t_{\max} \approx 0.
    \label{(Eq.B4)}
\end{equation}

Which, using the logarithm in \eqref{(Eq.B4)}, is equivalent to the
expression

\begin{align}
    \ln\left( k_{w} \right) & = \ln\left( t_{\max} \right)h + \ln\left( {\widetilde{k}}_{ap}\left( t = t_{\max} \right) \right),\label{(Eq.B5)}\\
    \overset{\Leftrightarrow}{}h & = \frac{1}{\ln\left( t_{\max} \right)}\ln\left( k_{w} \right) - \frac{\ln\left( {\widetilde{k}}_{ap}\left( t = t_{\max} \right) \right)}{\ln\left( t_{\max} \right)}.\label{(Eq.B6)}
\end{align}

\begin{figure}
\centering
\includegraphics[width=\linewidth]{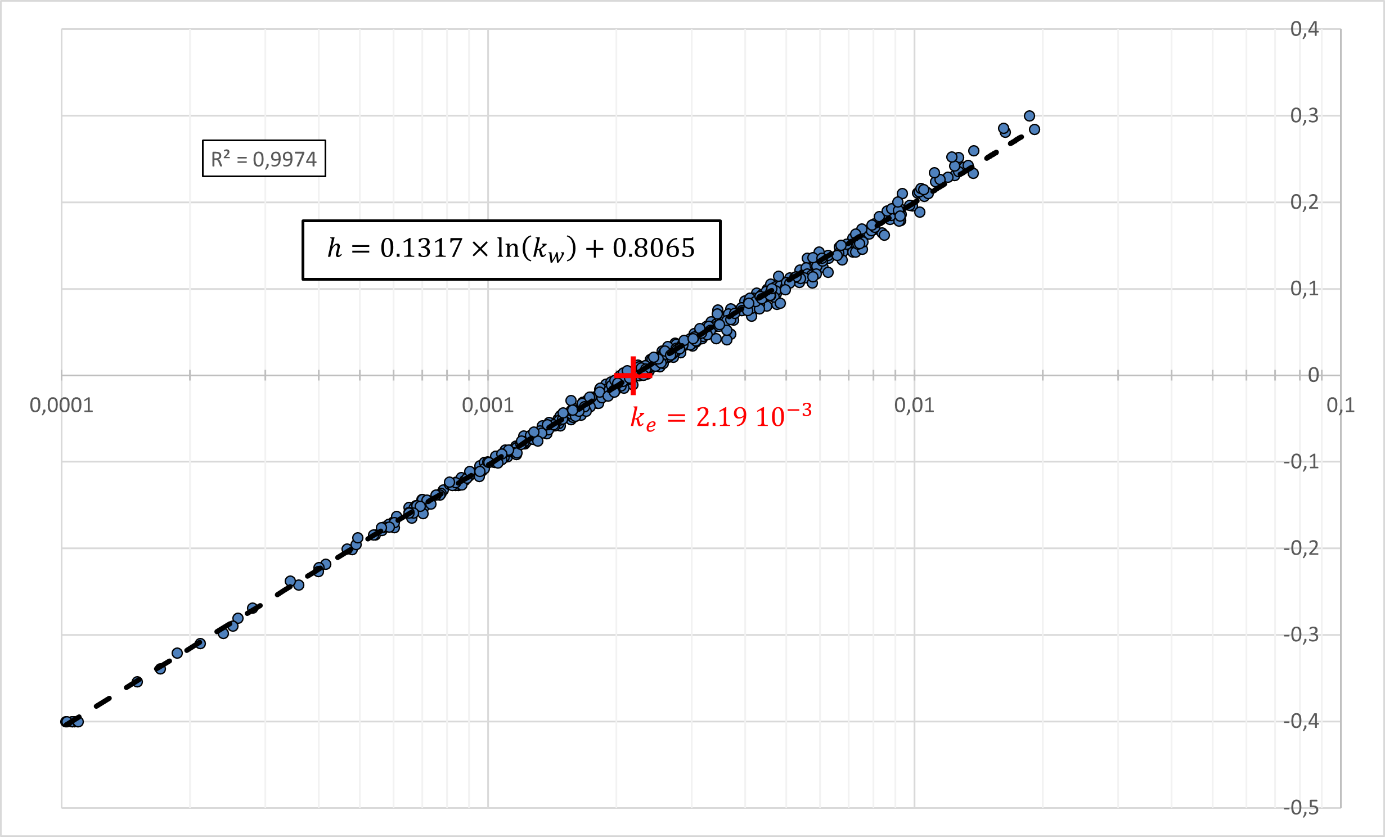}
\caption{Plot of the values of \(k_{w}\) and \(h\) (blue
points) in the Weibull model with the fitting described in this section
for the samples where the Weibull fitting is not perfect. The
trendline (black dashed line), its equation and the R² score associated
are given. The limit of the instantaneous rate coefficient \(k_{e}\)
(red mark) corresponding to \(k_{w}(h = 0)\) is displayed.}
\label{Fig.7}
\end{figure}

\eqref{(Eq.B6)} is exactly the relation in \eqref{(Eq.B1)} with
\(a = \frac{1}{\ln\left( t_{\max} \right)}\) and
\(b = - \frac{\ln\left( {\widetilde{k}}_{ap}\left( t = t_{\max} \right) \right)}{\ln\left( t_{\max} \right)}\).
If the instantaneous rate coefficient would not tend to the value
\(k_{e}\), since \(t_{\max}\) is identical in all measurements, we would
not find the logarithm relationship between \(h\) and \(k_{w}\)
displayed in Fig.\ref{Fig.7} as the intercept term would not be constant.
Thus\(\ {\widetilde{k}}_{ap}\left( t = t_{\max} \right)\) has to be a
constant that we define to be the value \(k_{e}\). With that definition,
\(k_{e}\) is the limit of the instantaneous rate coefficient
\(k_{ap}(t)\) that follows the general hydrolysis dynamic equation

\begin{equation}
    X_{t} = X_{\infty}\left( 1 - e^{- k_{ap} \times t} \right).\label{(Eq.B7)}
\end{equation}

Plotting \(h\) against \(\ln\left( k_{w} \right)\) and fitting Equation
(A1) onto it enables to find that \(k_{e}\) value, graphically being the
value of \(k_{w}\) where \(h = 0\), or by computing
\(k_{e} = e^{- \frac{b}{a}}\). For the MAGIC samples, we get
\(k_{e} = 2.19{\ 10}^{- 3}\ \min^{- 1}\) . The final experimental time
acquisition, \(t_{\max}\), can also be extracted from Fig.\ref{Fig.7} from the
slope of the curve. Here, we get an approximation of \(t_{\max}\) of
1987 min, compared to 1800 min that gives a theoretical slope of
\(\frac{1}{ln(1800)} = 0.1334\).

\section{Confidence intervals of \(k_{e}\) with the Weibull estimation}
\label{Appendix.C}
\setcounter{equation}{0}
\renewcommand{\theequation}{C.\arabic{equation}}
Based on \eqref{Eq.5} and \eqref{Eq.9}, we gave an estimation of \(k_{e}\), the
inherent reaction rate constant that reaches the instantaneous reaction
rate coefficient \(k_{a}(t)\) at \(t \rightarrow \infty\). Indeed, the
linear regression between the Weibull parameters \(h\) and
\(\ln\left( k_{w} \right)\) defined as

\begin{equation}
    h = a \times \ln\left( k_{w} \right) + b,\ \ a,b\ constants,
    \label{(Eq.C1)}
\end{equation}

Defines the estimator

\begin{equation}
    \widetilde{ke} = k_{w}(h = 0) = k_{w}^{0} = e^{\frac{- b}{a}}.
    \label{(Eq.C2)}
\end{equation}

For the rest of this section, we denote \(\ln\left( k_{w} \right)\) by
\(x\), i.e., \eqref{(Eq.C1)} \(\overset{}{\Leftrightarrow}h = a \times x + b\).

Let's consider a dataset of Weibull fits containing \(n\) couples of
\(h\) and \(x\) values:
\(\left\{ \left( h_{i},x_{i} \right)_{i \leq n} \right\}\). Coefficients
\(a\) and \(b\) are given by \(\widehat{a} = \frac{S_{xh}}{S_{xx}}\) and
\(\widehat{b} = \overline{h} - \widehat{a}\overline{x}\), with
\(S_{xx} = \sum_{}^{}\left( x_{i} - \overline{x} \right)^{2}\) ,
\(S_{xh} = \sum_{}^{}{(x_{i} - \overline{x})(h_{i} - \overline{h})}\).
The residual variance being
\({\widehat{\sigma}}^{2} = \frac{1}{n - 2}\sum_{}^{}\left( h_{i} - \overline{h} \right)^{2}\).

Finding a confidence intervals for \(\frac{- \widehat{b}}{\widehat{a}}\)
is equivalent of finding one for \(\widetilde{ke}\). Rather than trying
to define the distribution of \(z = \frac{- \widehat{b}}{\widehat{a}}\),
which is the ratio of two correlated random variables, we look for
values of \(z\) such that the linear regression evaluated in \(z\) is
compatible with \(h = 0\).

Observing that
\(\mathbb{V}\left( \widehat{a}z + \widehat{b} \right) = {\widehat{\sigma}}^{2}(\frac{1}{n} + \frac{\left( z - \overline{x} \right)^{2}}{S_{xx}})\),
once get that
\(T(z) = \frac{\widehat{a}z + \widehat{b}}{{\widehat{\sigma}}^{2}(\frac{1}{n} + \frac{\left( z - \overline{x} \right)^{2}}{S_{xx}})}\sim t_{n - 2}\),
i.e., the Student distribution with \(n - 2\) degrees of freedom. Then a
confidence interval, of level \(1 - \alpha\), for \(z\) is defined by
\(\left\{ z:\ \left| T(z) \right| \leq t_{n - 2,\frac{\alpha}{2}} = t^{*} \right\}\).
Solving \(T(z)^{2} = {t^{*}}^{2}\), and denoting
\(g = \frac{{t^{*}}^{2}{\widehat{\sigma}}^{2}}{{\widehat{a}}^{2}S_{xx}}\),
we get the confidence interval
\(\left\lbrack z_{1},z_{2} \right\rbrack\), where

\begin{equation}
    z_{1,2} = - \frac{\widehat{b}}{\widehat{a}} \times \frac{1}{1 - g} \pm \frac{t^{*}\widehat{\sigma}}{\widehat{a}(1 - g)}\sqrt{\frac{1}{n} + \frac{\frac{{\widehat{b}}^{2}}{{\widehat{a}}^{2}}}{S_{xx}} - g \times \frac{1}{n}}.
    \label{(Eq.C3)}
\end{equation}

\section{Examples of fitting for one sample of WWS and HAWS}
\label{Appendix.D}
\begin{figure}
\centering
\includegraphics[width=\linewidth]{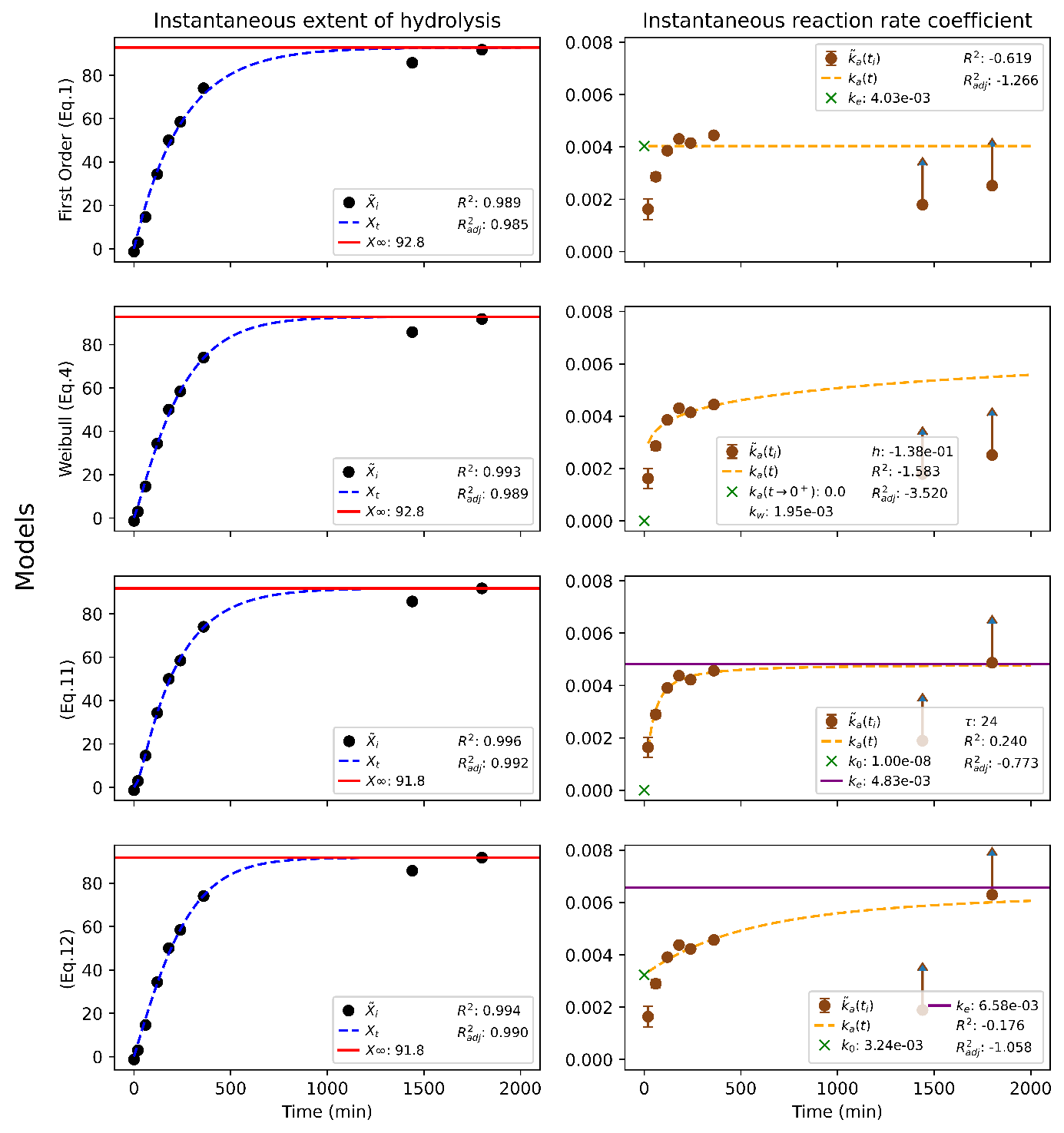}
\caption{Example of models\textquotesingle{} fit with data
from a WWS sample. Each row is a different model. The left panels detail
the calibration based on the hydrolysis data (\(X_{t}\)), while the
right panels detail the calibration based on the instantaneous reaction
rate associated (\(k_{a}(t)\)). Confidence intervals with a single, arrow-shaped 
ascending branch arise when the their length is higher than $3.10^{-3}$ and are kept
at this size to prevent the plot from becoming flattened.}
\label{Fig.8}
\end{figure}

The WWS sample in Fig.\ref{Fig.8} exhibits poor results with \(R^{2}\) score of
-0.619 (\(R_{adj}^{2}\): -1.266), -1.583 (\(R_{adj}^{2}\): -3.520),
0.240 (\(R_{adj}^{2}\): -0.773), and -0.176 (\(R_{adj}^{2}\): -1.058)
for \eqref{Eq.1}, \eqref{Eq.4}, \eqref{Eq.11}, and \eqref{Eq.12} respectively. Only the model
from \eqref{Eq.11} has a positive \(R^{2}\) score that is not satisfying
nonetheless. These results were to be expected, since it is clear from
Fig.\ref{Fig.1}c that the instantaneous reaction rate coefficient is not monotonic
in time. For all WWS samples \(k_{a}\) is increasing for around 300
minutes, and then decreasing over time to an asymptotic value of around
2.41 \(10^{- 3}\). In that regard, since none of the four models
presented here can simulate a non-monotonic dynamic for \(k_{a}\), none
should yield satisfying results, thus a new model for the instantaneous
reaction rate coefficient must be designed to fit the WWS example.

\begin{figure}
\centering
\includegraphics[width=\linewidth]{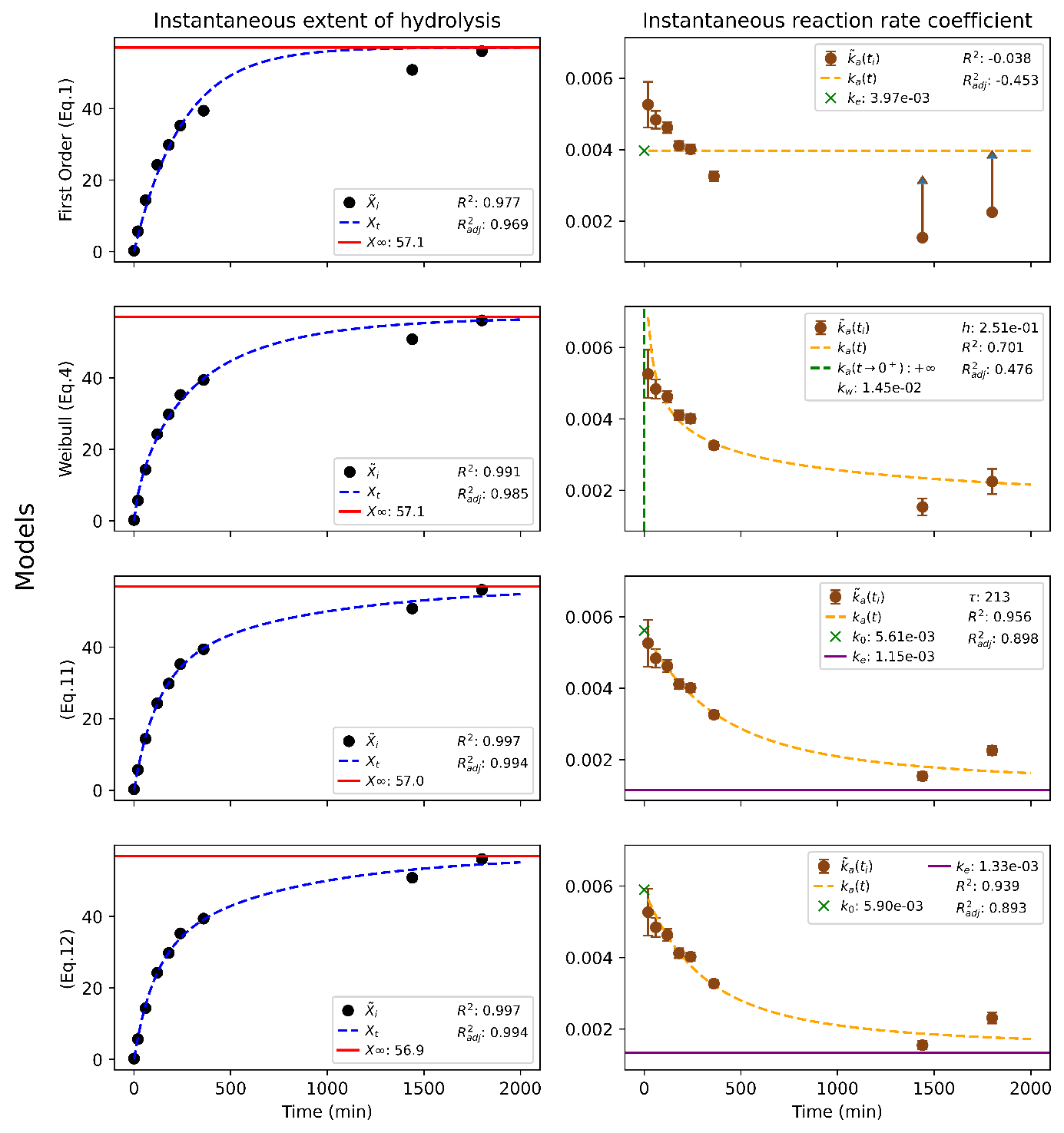}
\caption{Example of models\textquotesingle{} fit with data
from a HAWS sample. Each row is a different model. The left panels
detail the calibration based on the hydrolysis data (\(X_{t}\)), while
the right panels detail the calibration based on the instantaneous
reaction rate associated (\(k_{a}(t)\)). Confidence intervals with a single, arrow-shaped 
ascending branch arise when the their length is higher than $3.10^{-3}$ and are kept
at this size to prevent the plot from becoming flattened.}
\label{Fig.9}
\end{figure}

The HAWS sample in Fig.\ref{Fig.9} exhibits similar results as the ones observed
in Fig.\ref{Fig.5}, with the Weibull model having a higher \(R^{2}\) score of
0.701 (\(R_{adj}^{2}\): 0.476), but \eqref{Eq.11} and \eqref{Eq.12} have an even
higher \(R^{2}\) score of 0.956 (\(R_{adj}^{2}\): 0.898) and 0.939
(\(R_{adj}^{2}\): 0.893) respectively. Based on those
\(R_{adj}^{2}\ \)scores, \eqref{Eq.11} and \eqref{Eq.12} are equivalent to model the
HAWS sample example.

\printcredits

\bibliographystyle{cas-model2-names}

\bibliography{references}



\end{document}